\documentclass[12pt, reqno]{amsart}
\usepackage{times}
\usepackage{amssymb}
\usepackage{graphicx}
\usepackage{amsmath,amstext}
\newtheorem{thm}{Theorem}[section]
\newtheorem{cor}[thm]{Corollary}
\newtheorem{lem}[thm]{Lemma}

\numberwithin{equation}{section}
\newtheorem{conjecture}{\sffamily\bfseries Conjecture}
\date{}

\newcommand{\ip}[1]{\langle #1 \rangle}

\def\D{\mathbb{D}}

\begin{document}

\title[Commutant and Reducing subspaces  ]{Multiplication operators on the Bergman space via analytic continuation}

\author[Douglas]{Ronald G. Douglas}
\address{Department of Mathematics, Texas A$\&$M University, College Station, Texas 77843}
\email{rdouglas@math.tamu.edu}
\author[Sun]{Shunhua Sun}
\address{Institute of Mathematics, Jiaxing University, Jiaxing, Zhejiang, 314001, P. R. China}
\email{shsun$_{-}$cn@yahoo.com.cn}
\author[Zheng]{Dechao Zheng}
\address{Department of Mathematics, Vanderbilt University, Nashville, TN 37240}
\email{dechao.zheng@vanderbilt.edu}
\thanks{}

\begin{abstract} In this paper, using the group-like property of local inverses of a finite
Blaschke product $\phi$, we will show that the largest $C^*$-algebra in  the commutant of the
multiplication operator $M_{\phi}$ by $\phi$ on the Bergman space
 is finite dimensional, and its dimension equals the
number
  of connected components of the Riemann surface of $\phi^{-1}\circ\phi $ over the unit disk.
If the order of the Blaschke product $\phi$ is less than or equal to eight, then every $C^*$-algebra contained in
the commutant of $M_{\phi}$ is
abelian and hence the number of minimal reducing subspaces of $M_{\phi}$ equals the
number
  of connected components of the Riemann surface of $\phi^{-1}\circ\phi $ over the unit disk.

\end{abstract}


\maketitle

\section{Introduction}

Let $ \D$ be the open unit disk in $\mathbb C$. Let $dA$
denote Lebesgue area measure on the unit disk $\D$, normalized so
that the measure of~$ \D$ equals $1$. The Bergman space $L^{2}_{a}$
is the Hilbert space consisting of the analytic functions on~$\D$
that are also in the space $L^{2}( \D, dA)$ of square integrable
functions on $ \D$. For a bounded
analytic function $\phi$ on the unit disk, the multiplication
operator $M_{\phi}$ is defined on the Bergman space $L^{2}_{a}$
given by
$$M_{\phi}h=\phi h$$
for $h\in L^{2}_{a}.$

The classification of invariant subspaces or reducing subspaces of
various operators acting on function spaces has proved to be one
very rewarding research problem in analysis. Not only     has the
problem itself turned out to be important, but also the methods used
to solve it are interesting. The classical Beurling theorem
\cite{Beu} gives a complete characterization of the invariant
subspaces of the unilateral shift.  Extensions of this idea  have
led to many important works by other investigators. On the Bergman
space, the lattice of invariant subspaces of the Bergman shift
$M_{z}$ is huge and rich \cite{BFP} although a Beurling-type theorem
is established in \cite{ARS}.

A reducing subspace $\mathcal M$ for an operator $T$ on a Hilbert
space $H$ is a subspace $\mathcal M$ of $H$  such that $T{\mathcal
M}\subset {\mathcal M}$ and $T^{*}{\mathcal M} \subset {\mathcal
M}$. A reducing subspace $ {\mathcal M}$ of $T$ is called minimal if
the only reducing subspaces contained in $ {\mathcal M}$ are
${\mathcal M} $ and $\{0\}.$  Let $\{M_{\phi}\}^{\prime}$ denote the
commutant of $M_{\phi}$, which is the set of bounded operators on
the Bergman space commuting with $M_{\phi}.$ The problem of
determining the reducing subspaces of an
 operator is equivalent to finding the projections in the commutant of
 the operator. An $n$th-order Blaschke product is the function analytic on the unit disk $\D$ given by
$$\phi (z)=\prod_{j=1}^{n}\frac{z-a_{j}}{1-\overline{a}_{j}z}$$
for $a_{j}\in \D$.
  For an $n$th-order Blaschke product $\phi$,
 since $M_{\phi}^{*}$ belongs to the Cowen-Douglas class $B_n(\D),$
 \cite{CoD}, \cite{Do2},
one can apply   results from complex geometry to note that reducing
subspaces correspond to some subspaces of a single fiber which is
isomorphic to ${\mathbb C}^n.$  This implies immediately that there
can't  be more than $n$ pairwise orthogonal reducing subspaces of
$M_\phi$. However, the lattice of reducing subspaces of the bounded 
operator $M_\phi$ could still be infinite

Let ${\mathcal A}_{\phi}$ be the von Neumann algebra,  defined to be
the intersection of the commutants $\{M_{\phi}\}^{\prime}$ and
$\{M_{\phi}^{*}\}^{\prime}$. The goal of this paper is to study
${\mathcal A}_{\phi}$ for a Blaschke product $\phi$ of finite
order. This is a continuation of the investigation begun in
\cite{SZZ1}, \cite{SZZ2}. In \cite{SZZ1}, \cite{SZZ2}, one used the
Hardy space of the bidisk to study  multiplication operators on
$L_a^2 $ by bounded analytic functions on the unit disk $ \D$ and to
give complete classification of the reducing subspaces of
multiplication operators on $L_a^2 $ by Blaschke products with order
$3$ or $4$.

In this paper, we will take a completely different approach from the one in \cite{SZZ1}, \cite{SZZ2}.
On one hand,   the multiplication operators have deep connections with  the analytic properties of their symbols $\phi$.
 On the other hand, those multiplication operators  are typical subnormal
 operators whose minimal normal extensions have a thick spectrum  and the adjoints of the multiplication operators
 are in the Cowen-Douglas class \cite{CoD}.  We will make use of two more ingredients. One is local inverses of a finite
 Blaschke product and their analytic continuations on a subset of the unit disk. The germs induced by these local
 inverses have a group-like property by compositions of germs (for details of compositions of germs, see \cite{Vec}).
The group-like property was  used in \cite{Co1}, \cite{Th1},
 \cite{Th2} in studying the commutant of Toeplitz operators on the Hardy space. 
 
 The other
ingredient is the theory of subnormal operators \cite{Con} which, combined with
properties of the Cowen-Douglas classes \cite{CoD}, can be used to show 
   that unitary operators
 in the commutant of the multiplication operators have a nice representation. Combining the group-like property of local
 inverses and
 the nice representation of unitary operators, we will obtain a symmetric and unitary matrix representation of
 the action of these unitary operators acting on reproducing kernels.

 Our main result in the paper is that the dimension of  ${\mathcal A}_{\phi}$    equals  the
number  of connected components of the Riemann surface of
$\phi^{-1}\circ\phi $ over $ \D$ .
  This result was obtained for Blaschke products of order $3$ or $4$ in \cite{SZZ1},
\cite{SZZ2} and   suggests the following conjecture.
\begin{conjecture}\label{con1}
 For a Blaschke product $\phi$ of finite order,
 the number of nontrivial minimal reducing subspaces of $M_\phi$   equals  the number
  of connected components of the Riemann surface of $\phi^{-1}\circ\phi $ over $\D$.
  \end{conjecture}

 Our main result implies that Conjecture \ref{con1} is equivalent to that the $C^*$-algebra ${\mathcal A}_{\phi}$ is abelian. The conjecture is also equivalent to the statement whether or not the minimal reducing subspaces are orthogonal.
 For a Blaschke product $\phi$ with order smaller than or equal to $8$, we will confirm the conjecture by showing that ${\mathcal A}_{\phi}$ is abelian in the last section.


This paper is motivated  by Richter's work on unitary
equivalence of invariant subspaces of the Bergman space
  \cite{Ric}, Stephenson's work on hypergroups  \cite{St1}, \cite{St2}, \cite{St3} and Zhu's conjecture
  on the number of minimal reducing subspaces of
  $M_{\phi}$ \cite{Zh1}. Many ideas in the paper are, however, inspired by nice works on the commutant of
 an analytic Toeplitz operator  on the Hardy space  in \cite{Co1}, \cite{Th1},
 \cite{Th2}. C. Cowen also used the Riemann surface of
$\phi^{-1}\circ\phi $ over $ \D$ to describe the commutant of the multiplication operator by $\phi$ on the Hardy space
 in \cite{Co1}.

  Thomson's  representation of local inverses in \cite{Th1}, \cite{Th2} is also useful in the Bergman space context.
Using those local inverses  as  in \cite{Th1}, one can easily see
that for an analytic and nonconstant function $f$ in the closed unit
disk $\overline{\D}$, there exists a finite Blaschke product $\phi$
  such that
  $$\{M_{f}\}^{\prime}=\{M_{\phi}\}^{\prime}.$$
This shows that the multiplication operator by a finite Blaschke
product will play an important
  role in studying the other multiplication operators on the Bergman space.

We would like to point out that the results and arguments can carry
over to the weighted Bergman spaces, but do not work on the Hardy
space.   On the Hardy space,  because the spectral measure of the
minimal normal extension of the multiplication operator by $\phi$ is
supported on the unit circle which is its essential spectrum. On the
other hand, the spectral measure of the minimal normal extension of
the multiplication operator by $\phi$ on the weighted Bergman space
is supported on the unit disk, which is its spectrum.  
The key fact is that the
spectral measure is supported on the boundary of the disk in view of the
maximum principle.  Although the proof doesn't make it explicit, we
believe if points on the interior of the disk are essential with respect
to the spectral measure then the argument goes through.

\section{Analytic continuation and Local inverses}

First we introduce some notation. An analytic function element is a
pair $(f, U)$, which consists of  an open disk $U$ and an analytic
function $f$ defined on this disk. A finite sequence ${\mathcal
U}=\{(f_{j}, U_{j})\}_{j=1}^{m}$ is a continuation sequence if
\begin{itemize}

\item $U_{j}\cap U_{j+1}$
is not empty for $j=1, \cdots , m-1$ and

\item   $f_{j}\equiv f_{j+1}$ on $U_{j}\cap U_{j+1}$, for $j=1, \cdots , m-1$.

\end{itemize}
Let $\gamma$ be an arc with parametrization $z(t)$,  $z(t)$ being a continuous function on an interval $[a, b]$. A sequence $\{U_{1}, \cdots , U_{m}\}$ is admissible or a covering chain for $\gamma $ if each $U_{j}$ is an open disk, and if there exist
increasing numbers $t_{1}, \cdots t_{m}$ in $[a, b]$ such that $z(t_{j})\in U_{j}$ for $j=1, \cdots , m$ and
$$z(t)\in \left\{\begin{array}{ll}
              U_{1},&  ~a\leq t\leq t_{1}\\
              U_{j}\cup U_{j+1}, &   t_{j}\leq t \leq t_{j+1}\\
              U_{m}, &  t_{m}\leq t\leq b.
              \end{array} \right.
              $$

              A continuation sequence ${\mathcal U}=\{(f_{j}, U_{j})\}_{j=1}^{m}$ is an analytic continuation
               along the arc $\gamma$ if the sequence $U_1 , \cdots , U_{m}$ is admissible for $\gamma .$ Each of
               the elements $\{(f_{j}, U_{j})\}_{j=1}^{m}$
        is an analytic continuation of the other along the curve $\gamma$. We say that the analytic function
        $f_{1}$ on     $U_{1}$ admits a continuation to $U_{m}.$ A
        famous result on analytic continuations is the following Riemann
        monodromy theorem \cite{Ahl},  \cite{Gri} and \cite{Vec}.

\begin{thm}\label{RMT} Suppose $\Omega \subset {\mathbb C}$ is a simply
connected open set. If an analytic element, $(f, U)$ can be
analytically continued along any path inside $\Omega$, then this
analytic function element can be extended to be a single-valued
holomorphic function defined on the whole of $~~\Omega .$
\end{thm}

Let $\phi$ be an $n$-th order  Blaschke product. Let
$$E=\D\backslash [\phi^{-1}(\phi (\{\beta\in \D:\phi^{\prime}(\beta )=0\}))].$$
Note that $\D\backslash$ is finite. For an open set $V\subset \D$,
we define a local inverse of $\phi$ in $V$ to be a function $f$
analytic in $V$ with $f(V)\subset \D$ such that $\phi (f(z))= \phi
(z)$ for every $z$ in $V$.
 That is, $f$ is  a branch of $\phi^{-1}\circ\phi $ defined in $V$.

 A finite collection, $\{f_{i}\}$,  of local inverses in $V$ is complete if for each $z$ in $V$,
$$\phi^{-1}(\phi (z))\cap \D\subset \{f_{i}(z)\}$$
and
$$f_{i}(z)\neq f_{j}(z)$$
for $i\neq j$. An open set $V$ is invertible if there exists a complete collection of local inverses in $V$.

A local inverse $ (f, V)$ admits an analytic continuation along the
curve $\gamma$ in $E$ if there is a continuation sequence ${\mathcal
U}=\{(f_{j}, U_{j})\}_{j=1}^{m}$ admissible for $\gamma$ and
$(f_{1}, U_{1})$ equals $(f, V).$  A local inverse in $V\subset E$
is admissible for $\phi$ if it admits unrestricted continuation in
$E$. Note that the identity function is always admissible. The set
of admissible local inverses has the useful property that it is
closed under composition, which can be shown as follows. Let $f$ and
$g$ be admissible local inverses in open discs $V$ and $W$ centered
at $a$ and $f(a)$, respectively, with $f(V)\subset W$. Let $\gamma$
be a curve in $E$ with initial point $a$. Since $f$ is admissible,
it can be analytically continued along $\gamma$. There is an obvious
image curve $\tilde{\gamma}$ of $\gamma$ under this analytic
continuation along $\gamma$. Since $g$ is also admissible, it can be
analytically continued along $\tilde{\gamma}$. By refining the
covering chain of $\gamma$, if necessary, we can assume that if
$\tilde{V}$ is a covering disc of $\gamma$ and $ (\tilde{f},
\tilde{V})$ the corresponding function element, then
$\tilde{f}(\tilde{V})$ is contained in one of the covering discs of
$\tilde{\gamma}$. We now compose corresponding function elements in
the analytic continuations along $\gamma$ and $\tilde{\gamma}$ to
obtain an analytic continuation for $(g\circ f, V)$ along $\gamma$.

  Let $V$ be an invertible open disc and let $\{f_{i}\}$ be the family of admissible local inverses in $V$.
  By shrinking $V$, we can assume that
  each $f_{i}(V)$ is contained in an invertible open disc $W_{i}$. Let $\{g_{ij}\}_{j}$ be the family of admissible
  local inverse in $W_{i}.$ Since $g_{ij}\circ f_{i}$ is admissible and $g_{ij}\circ f_{i}\neq g_{ik}\circ f_{i}$
  if $j\neq k$, we observe that $\{g_{ij}\circ f_{i}\}=\{f_{j}\}$, for each $i$. In particular,
  for each $f_{i}$, there exists $g_{ij}$ such that $g_{ij}\circ f_{i}$ is the identity function in
  $V$, which means $f_{i}^{-1}=g_{ij}$ for some $j$,
  and thus $f_{i}^{-1}$
  is admissible.

   For each $z\in E$, the function $\phi$ is one-to-one in some open neighborhood $D_{z_{i}}$ of each point $z_{i}$ in  $\phi^{-1}\circ\phi (z)=\{z_{1}, \cdots , z_{n}\}$.  Let  $\phi^{-1}\circ\phi=\{\rho_{k}(z)\}_{k=1}^{n}$ be $n$ solutions $\phi (\rho (z))=\phi (z).$ Then $\rho_{j}(z)$ is
 locally analytic and arbitrarily continuable in $E$. Assume that $\rho_{1}(z)=z$. Every open subset $V$ of $E$ is
 invertible for $\phi$.
 Then $\{\rho_{j}\}_{j=1}^{n}$ is the family of  admissible local inverses in some invertible open disc $V\subset \D$.
 For a given point $z_{0}\in V$, label those local inverses as $\{\rho_{j}(z)\}_{i=1}^{n}$ on $V$.
 If there is a loop $\gamma$ in $ E$ at $z_{0}$ such that $\rho_{j}$ and  $\rho_{j^{\prime}}$ in $\{\rho_{i}(z)\}_{i=1}^{n}$
 are mutually analytically continuable along $\gamma$, we can then write
$$\rho_{j}\thicksim \rho_{j^{\prime}},$$
and it is easy to check that $\thicksim$ is an equivalence relation.
Using this equivalence relation, we partition
$\{\rho_{i}(z)\}_{i=1}^{n}$ into equivalence classes
$$\{G_{i_1}, G_{i_{2}}, \cdots , G_{i_{q}}\},$$ where
$i_{1}=1<i_{2}<i_{3}<\cdots <i_{q}\leq n$
 for some integer $1<q\leq n$ and $\rho_{i_{k}}$ is in $G_{i_{k}}.$
Not all
 of the branches of $\phi^{-1}\circ\phi$ can be continued to a different branch. For
 example,
 $z$ is a single valued branch of $\phi^{-1}\circ\phi$. Then $q$ is
 greater than $1$.
 Thus each element in $G_{i_{k}}$ extends analytically to the other element in  $G_{i_{k}}$, but it does not extend to
 any element in
 $G_{i_{l}}$ if $i_{k}\neq i_{l}$.
So
 \begin{equation}\label{Concom}
 \{\rho_{i}(z)\}_{i=1}^{n}=\cup_{k=1}^{q}G_{i_{k}}.
 \end{equation}
The  collection $\{G_{i_1}, G_{i_{2}}, \cdots , G_{i_{q}}\}$ does
not depend on the choice of $z_{0}$ in $E$.

Let $k_{\alpha}$ denote the reproducing kernel of the Bergman space
at the point $\alpha$ in $\D$.  As in \cite{DoM}, \cite{Th1}, we will
use local inverses to obtain a local representation of an operator
$T$ in the commutant $\{M_{\phi}\}^{\prime}$. The proof of the
following theorem is similar to the ones in \cite{Co1}, \cite{Th1}.

  \begin{thm}\label{localrepu} Let $\phi$ be a finite Blaschke product,  $U$ be an invertible
  open set of $E$, and t $\{\rho_{i}(z)\}_{j=1}^{n}$  be a complete collection of local inverses  on $U$. Then for  each $T$ in $\{M_{\phi}\}^{\prime}$, there are analytic functions $\{s_{i}(\alpha )\}_{i=1}^{n}$
  on $U$ such that for each $h$ in the Bergman space $L^{2}_{a},$
  \begin{eqnarray*}
  Th(\alpha )&= &\sum_{i=1}^{n}s_{i}(\alpha )h(\rho_{i}(\alpha )),\\
  T^{*}k_{\alpha}&=&\sum_{i=1}^{n}\overline{s_{i}(\alpha )}k_{\rho_{i}(\alpha )}
  \end{eqnarray*}
 for each $\alpha$ in $U$. Moreover, these functions $\{s_{i}(\alpha )\}_{i=1}^{n}$ admit  unrestricted continuation in $E$.
  \end{thm}

  \proof Since $T$ commutes with $M_{\phi}$, the adjoint $T^{*}$ commutes with $M_{\phi}^{*}$. Thus $T^{*}$ commutes with $M_{\phi -\phi (\alpha )}^{*}$
  for each $\alpha$ in the invertible set $V$. So the kernel of $M_{\phi -\phi (\alpha )}^{*}$ is invariant for $T^{*}.$ Note that the kernel of $M_{\phi -\phi (\alpha )}^{*}$ is the finite dimensional space spanned by $\{k_{\rho_{i}(\alpha )}\}_{i=1}^{n}$. Hence for each $\alpha$ in $U$,
  there is a sequence $\{s_{i}(\alpha )\}_{i=1}^{n}$ of complex numbers such that
  $$T^{*}k_{\alpha}=\sum_{i=1}^{n}\overline{s_{i}(\alpha )}k_{\rho_{i}(\alpha )}.$$
  Thus for each $h$ in $L^{2}_{a},$ we have
  \begin{eqnarray*}
  Th(\alpha ) &=&\ip{Th, k_{\alpha}}\\
              &=&\ip{h, T^{*}k_{\alpha}}\\
              &=&\ip{h, \sum_{i=1}^{n}\overline{s_{i}(\alpha )}k_{\rho_{i}(\alpha )}}\\
              &=& \sum_{i=1}^{n}s_{i}(\alpha )\ip{h, k_{\rho_{i}(\alpha )}}\\
              &=& \sum_{i=1}^{n}s_{i}(\alpha )h(\rho_{i}(\alpha )).
  \end{eqnarray*}
  To finish the proof we need to show that $\{s_{i}(z)\}_{i=1}^{n}$ are analytic in $U$. To do so, for each $i$, define
  $$P_{i}(\alpha , z)=\prod_{j\neq i}^{n}(z-\rho_{j}(\alpha ))$$
  for $z$ in $\D$ and $\alpha$ in $U$. Thus $\{P_{i}(\alpha , z)\} $ is a family of functions  analytic in $z$ on $\D$ and
  analytic in $\alpha$  on $U$.
  An easy calculation gives that for each $\alpha$ in $U$
  \begin{eqnarray*}
  \ip{P_{i}(\alpha , .), T^{*}k_{\alpha}}&=& \ip{P_{i}(\alpha , .), \sum_{j=1}^{n}\overline{s_{j}(\alpha )}k_{\rho_{i}(\alpha )}}\\
                                         &=& \sum_{j=1}^{n}s_{j}(\alpha )\ip{P_{i}(\alpha , .),k_{\rho_{i}(\alpha )}}\\
                                         &=& s_{i}(\alpha )P_{i}(\alpha , \rho_{i}(\alpha ))\\
                                         &=& \prod_{j\neq i}^{n} (\rho_{i}(\alpha )-\rho_{j}(\alpha )).
                                         \end{eqnarray*}
  Thus
  $$s_{i}(\alpha )=\frac{\ip{P_{i}(\alpha , .), T^{*}k_{\alpha}}}{\prod_{j\neq i}^{n} (\rho_{i}(\alpha )-\rho_{j}(\alpha ))}$$
  for $\alpha$ in $U$ and hence  $s_{i}$ is analytic in $U$. Noting that $\{\rho_{i}(z)\}_{j=1}^{n}$ admit  unrestricted continuation in $E$, we conclude that
  the functions $\{s_{i}(\alpha )\}_{i=1}^{n}$ admit  unrestricted continuation in $E$ to   complete  the proof.

 \section{\label{Riems}Riemann surfaces $\phi^{-1}\circ\phi $ over $\D$}

Let $\phi=\frac{P(z)}{Q(z)}$ be an $n$-th order Blaschke product of
where $P(z)$ and $Q(z)$ are two coprime polynomials of degree less
than or equal to $n$. In this section we will study the Riemann
surface for the Blaschke product $\phi^{-1}\circ\phi$ over $\D$. In
particular, it was shown in \cite{SZZ2} how the number of the
connected components of the Riemann surface $\phi^{-1}\circ\phi$
over $\D$ is related to the zeros of $\phi$ for the fourth order
Blashcke product $\phi$. Let
$$f(w,z)=P(w)Q(z)-P(z)Q(w).$$
 Then $f(w,z)$ is a polynomial of $w$ with degree $n$ and its coefficients are polynomials
of $z$ with degree  $n$. For each $z\in \D$, $f(w,z)=0$ has exactly
$n$ solutions in $\D$ counting multiplicity. An algebraic function
is a function $w=g(z)$ defined for values $z$ in $\D$ by an equation
$f(w,z)=0.$

Let $ {\mathcal C}$ denote
 the set of the critical points of $\phi$ in $\D$ and
 $$ {\mathcal F}=\phi^{-1}\circ\phi ( {\mathcal C})=\{z_{1}, \cdots , z_{m}\}.$$
 Then $ {\mathcal F}$ is a finite set and is called the set of branch points of $\phi$, and
 $\phi^{-1}\circ\phi=\{\rho_{k}(z)\}_{k=1}^{n}$ is an $n$-branched analytic
 function defined and arbitrarily continuable in ${\D}\backslash {\mathcal F}$. Not all
 of the branches of $\phi^{-1}\circ\phi$ can be continued to a different branch, for example
 $z$ is a single valued branch of $\phi^{-1}\circ\phi$. The Riemann surface $S_{\phi}$ for
 $\phi^{-1}\circ\phi $ over $\D$ is an $n$-sheeted cover of $\D$ with at most $n(n-1)$ branch points, and
 it is not connected.
 We denote a point of $S_{\phi}$ lying over ${\D}\backslash{\mathcal F}$ by
 $(\rho (\alpha ), \alpha )$, where $\alpha$ is in ${\D}\backslash{\mathcal F}$ and $\rho$ is
 a branch of $\phi^{-1}\circ\phi$ defined in a neighborhood of
 $\alpha$.

 Visualization of Riemann surfaces is complicated by the fact that they are embedded in ${\mathbb C}^{2}$,
 a four-dimensional real space. One aid to constructing and visualizing them is a method
known as ``cut and paste''. Here we present only
  details on how to construct $S_{f}.$ For general cases, see \cite{Art}, \cite{Bli}, \cite{Fro}. We begin with $n$
  copies of the unit disk $\D$, called sheets. The sheets
 are labeled $\D_{1}, \cdots , \D_{n}$ and stacked up over $\D$. Then $\{z_{1}, \cdots , z_{m}\}$
 are the branch points. Suppose $\Gamma$ is a curve drawn through those branch points and a fixed point on the
 unit circle so that $\D\backslash\Gamma$
 is a simply connected region. By the Riemann monodromy theorem, $n$ distinct function elements $\rho_{k}(z)$,
 $k=1, \cdots , n$ of the algebraic equation
 $$f(w,z)=0$$
 can be extended to be a single-valued holomorphic functions defined over the whole of $\D\backslash\Gamma$. We denote these extended functions still
 by $\rho_{j}(z)$.  We may assume that $\Gamma$ consists of line segments
 $l_{k}$ to connect $z_{k}$ to $z_{k+1}.$ The sheets $\D_{j}$ are cut open along those line segments $l_{k}$.
 Then various sheets are glued to others along opposite edges of cuts. With the point in the $k$-th sheet over
 a value $z$ in $\D\backslash\Gamma$ we associate the pair of values $(\rho_{k}(z), z)$.
 In this way a one-to-one correspondence is set up between the points in $S_{f}$ over $\D\backslash\Gamma$
and the pair of points on the $n$ sheets over $\D\backslash\Gamma$.
 In order to make the correspondence continuous along the cuts exclusive of their ends, let two regions
 $R_{1}$ and $R_{2}$ be defined in a neighborhood of each cut $l_{i}$. On each of the $n$ sheets, in the region
 formed by $R_{1},$ $R_{2}$ and the cut $l_{i}$ between them exclusive of its ends, the values of the algebraic function
 $w=g(z)$ form again $n$ distinct holomorphic functions $\rho_{k}(z)$ ($k=1, \cdots , n$), and these can be   numbered
 so that $g_{l}(z)=\rho_{l}(z)$ in $R_{1}.$ In the region $R_{2}$ the functions $g_{k}(z)$
   are the same functions in the set  $\{\rho_{k}(z)\}$ but possibly in a different order. We join the edge of the cut bounding $R_{1}$ in the $k$-th
  sheet to the edge bounding $R_{2}$ in the $l$-th sheet, where $l$ is so determined that $g_{k}(z)=\rho_{l}(z)$
  in $R_{2}.$ The continuous Riemann surface so formed has the property that points in the Riemann surface $S_{f}$
  over non-branch points $\D\backslash\{z_{1}, \cdots , z_{m}\}$ are in one-to-one continuous correspondence with
  the nonsingular points $(w,z)$ which satisfies the equation
  $f(w,z)=0$. We not only get a
manifold; that is,   these identifications are continuous but the
Riemann surface also has an analytic structure or the match ups are
analytic.

   We need to use the number of connected components of the Riemann surface $S_{\phi}$ in the last two
   sections. By the unique factorization theorem for the  ring ${\mathbb C}[z,w]$  of
   polynomials in $z$ and $w$, we can factor
   $$f(w,z)=\prod_{j=1}^{q}p_{j}(w,z)^{n_{j}},$$
   where $p_{1}(w,z), \cdots , p_{q}(w,z)$ are irreducible
   polynomials.
   Bochner's Theorem \cite{Wal} says that $\phi$ has finitely many critical
  points in the unit disk $\D$. Thus we have
$$f(w,z)=\prod_{j=1}^{q}p_{j}(w,z).$$
   The following theorem
   implies that the number of connected components equals the number of irreducible factors $f(w,z).$ This result
   holds for Riemann surfaces over   complex plane (cf.  \cite{Bli}, page 78 and  \cite{Fro}, page
   374).

   \begin{thm}\label{Riem} Let $\phi (z)$ be an $n$-th order  Blaschke
   product and $f(w,z)=\prod_{j=1}^{q}p_{j}(w,z).$
   Suppose that  $p(w,z)$  is one of
   factors of $f(w,z)$.  Then the
   Riemann surface $S_{p}$ is connected if and only if $p(w,z)$ is
   irreducible. Hence $q$ equals the number of connected components of the Riemann surface
   $S_{\phi}=S_{f}$.
   \end{thm}

   \proof
   Let $\{z_{j}\}_{j=1}^{m}$ be the branch points of $p(w,z)=0$ in $\D$. Bochner's Theorem \cite{Wal} says that those
   points $\{z_{j}\}_{j=1}^{m}$ are contained in a compact subset of $\D.$
   Suppose that $p(w,z)$ is irreducible. If the
   Riemann surface $S_{p}$ is not connected, let $\{\rho_{k}(z)\}_{k=1}^{n_{p}}$ be $n_{p}$ distinct branches of $p(w,z)=0$
   over $\D\backslash\Gamma$.
   Then $\{\rho_{k}(z)\}_{k=1}^{n_{p}}$ are also roots of the equation
   $$\phi (w)-\phi (z)=0.$$

   Assuming that $S_{p}$ is not connected, we will derive a contradiction. Suppose that one connected component of $S_{p}$ is made up of the sheets corresponding to
   $\{\rho_{1}, \cdots , \rho_{n_{1}}\}$ ($n_{1}<n_{p})$. Let $\sigma_{s} (x_{1}, \cdots , x_{n_{1}})$
    be elementary symmetric functions of variables $x_{1},$ $ \cdots,$ $ x_{n_{1}}$ with degree $s$:
     $$\sigma_{s} (x_{1}, \cdots , x_{n_{1}})=\sum_{1\leq j_{1}<j_{2}<\cdots <j_{s}\leq n_{1}}x_{j_{1}}x_{j_{2}}\cdots x_{j_{s}}.$$
     Since the continuation of any path in ${\D}\backslash {\mathcal F}$ only leads to a permutation in
     $\{\rho_{1}(z), \cdots , \rho_{n_{1}}(z)\},$ every $\sigma_{s}(z)=\sigma_{s}(\rho_{1}(z), \cdots ,\rho_{n_{1}}(z))$
  is unchanged under such a permutation and hence  is a holomorphic function well-defined on
  $\D\backslash\{z_{j}\}_{j=1}^{m}$ and analytically extends on a neighborhood of
   the unit disk although $\rho_{j}(z)$ is defined only on $\D\backslash\Gamma$.

   Note that $\rho_{j}(z)$
   is in $\D$. Thus $\sigma_{s}(\rho_{1}(z), \cdots ,\rho_{n_{1}}(z))$ is bounded on $\D\backslash\{z_{j}\}_{j=1}^{m}$.
   By the Riemann removable singularities theorem,
   $\sigma_{s}(\rho_{1}(z), \cdots ,\rho_{n_{1}}(z))$ extends analytically on   $t\D$ for some $t>1$. Now we extend
   $\sigma_{s}(\rho_{1}(z), \cdots ,\rho_{n_{1}}(z))$ to the complex plane $\mathbb C$. For each $z\in {\mathbb C}\backslash\D$, define
   $$f_{s}(z)=\sigma_{s}(\frac{1}{\overline{\rho_{1}(\frac{1}{\bar{z}})}}, \cdots,\frac{1}{\overline{\rho_{n_{1}}(\frac{1}{\bar{z}})}}). $$
   By Theorem 11.1 on page 25 \cite{Bli}, near an ordinary point $z=a$,  each function $\rho_{j}(z)$ has a
   power series of  $z-a$. By Lemma 13.1 on page 29 \cite{Bli}, each function $\rho_{j}(z)$ has
   a Laurent series of a fractional power of $(z-a)$ but the number of terms with negative exponents must be finite.
   Thus $f_{s}(z)$ is a meromorphic function in ${\mathbb C}\backslash\D$ and
   the point at infinity is a pole of each $f_{s}(z)$.
   Moreover, $f_{s}(z)$ is  analytic in a neighborhood of the unit
   circle. So $f_{s}(z)$ is analytic in $t\D\backslash r\D$ for $0<r<1<t.$

   Next we show that  on the unit circle except for one point,
   $$\frac{1}{\overline{\rho_{i}(\frac{1}{\bar{z}})}}=\rho_{i}(z) $$
   for each $i$.
To do this, for each $\rho\in\{\rho_{i}\}_{i=1}^{n}$, noting that
$\phi$ is analytic on a neighborhood $\mathcal V$ of the closure of
the unit disk and $\rho$ extends analytically on the unit circle
minus $ \Gamma$, we have that for each $w\in {\mathcal V}$ with
$|z|=1$,
$$|\rho(z)|=1. $$
To prove the above fact, we follow an argument in \cite{Th1}. If
$|\rho (z)|<1$, as
$$\phi (\rho (z))=\phi (z)$$ and $z$ is not a branch point of $\phi$, then
$\rho (z)$ is in $\D\backslash\{z_{j}\}_{j=1}^{m}$. Thus $\rho^{-1}$
is a local inverse near $\rho (z)$ but its range is not contained in
the unit disk $\D$. So $\rho^{-1}$ is not admissible and neither is
$\rho$. Hence $|\rho (z)|=1.$ This means
$$\rho (z)=\frac{1}{\overline{\rho (\frac{1}{\bar{z}})}}$$
for $|z|=1.$

   Thus
   $$\sigma_{s}(\rho_{1}(z), \cdots ,\rho_{n_{1}}(z))=f_{s}(z).$$
   for $z$ on the unit circle except for one point. So
   $$\sigma_{s}(\rho_{1}(z), \cdots ,\rho_{n_{1}}(z))=f_{s}(z)$$
   in a neighborhood of the unit circle. Define
   $$F_{s}(z)=\{\begin{array}{cc}
              \sigma_{s}(\rho_{1}(z), \cdots ,\rho_{n_{1}}(z))&  ~z\in \bar{\D}\\
              f_{s}(z) &   ~z\in {\mathbb C}\backslash\D .
              \end{array}
              $$
          Thus $F_{s}(z)$ is a meromorphic function in ${\mathbb C}$ and the point at infinity
          is a pole of each $F_{s}(z)$. Hence $F_{s}(z)$ is a rational function of $z$ and
          so is   $\sigma_{s}(\rho_{1}(z), \cdots ,\rho_{n_{1}}(z))$ in $\D$.

        Now consider the polynomial
        $$f_{1}(w,z)=w^{n_{1}}-\sigma_{1}(z)w^{n_{1}-1}+\cdots +(-1)^{n_{1}}\sigma_{n_{1}}(z)=\prod_{j=1}^{n_{1}}(w-\rho_{j}(z))  $$
        whose coefficients are rational functions of $z$. Thus
        $$p(w,z)=f_{1}(w,z)f_{2}(w,z)$$
        for another polynomial $f_{2}(w,z)$. This implies that $p(w,z)$ is
        reducible, which  contradicts the assumption that $p(w,z)$ is irreducible.
        Hence the Riemann surface $S_{p}$ is connected.

        If $p(w,z)$ is reducible, noting that every root $\rho (z)$ of $p(w,z)$ is also a root of $\phi(w)-\phi (z)$,
         by  Bochner's Theorem \cite{Wal}, we see that  $p(w,z)$ does not have multiple roots for
         $z\in \D\backslash\{z_{j}\}_{j=1}^{m},$  we can factor  $$p(w,z)=p_{1}(w,z)\cdots p_{q_{1}}(w,z)$$
         for some irreducible polynomials $p_{l}(w,z)$ with degree
         $k_{l}.$

Let $\{\rho_{l1}(z), \cdots , \rho_{lk_{l}}(z)\}$ be roots of the
equation $p_{l}(w,z)=0$.  Thus
$$p_{l}(\rho_{lj}(z), z)=0.$$
From the identity theorem of analytic functions, each analytic
continuation $\hat{\rho}_{lj}(z)$ must still satisfy the equation
$$p_{l}(\hat{\rho}_{lj}(z), z)=0,$$
        and so $\hat{\rho}_{lj}(z)$ must be in $\{\rho_{lj}\}_{j=1}^{k_{l}}.$  Since $p_{l}$ is irreducible, by the above argument, we see that
          the continuations of the roots $\rho_{l1}(z),$
        $ \cdots ,$ $ \rho_{lk_{l}}(z)$ of
        $p_{l}(w,z)$ are always roots of $p_{l}(w,z)$. Hence crossing a cut permutes the set of roots, and the
        Riemann surface $S_{p_{l}}$ is connected. This gives that $S_{p}$ has $q_{1}$ connected components, one for each of factors $p_{1}, \cdots , p_{q_{1}}$.
       In particular, $S_{f}$ has $q$ connected components.  This completes the proof.

\section{\label{repUO} Representation of unitary operators}

In this section, we will obtain a better  representation of a
unitary operator in the commutant $\{M_{\phi}\}^{\prime}$ on
$E\backslash\Gamma$. It will be defined  using the orientation
$\{\{\rho_{i}\}_{i=1}^{n}, U\}$ of a complete collection of local
inverses $\{\rho_{i}\}$ on $U$ for a small invertible open set $U$
of  $E$. The theory of subnormal operators \cite{Con} plays also an
important role in this section.

 Next we need to order the set $\{\rho_{j}\}_{j=1}^{n}$ globally over a simply connected subset of $E$. To do this, take an invertible small open set $U$ of  $E$ such that the intersection of
  $\rho_{j}(U)$ and $ \rho_{k}(U)$ is empty for $j\neq k$. We can always do so by shrinking  $U$ sufficiently. In this section, we fix the small invertible  open set $U$.
  We label $\{\rho_{j}(z)\}_{j=1}^{n}$ as $$\{\rho_{1}(z), \rho_{2}(z), \cdots , \rho_{n}(z)\}$$ and assume $\rho_{1}(z)=z$.
  Take a curve $\Gamma$ through the finite set ${\mathcal F}$ and connecting a point on the unit circle so that
  $E\backslash\Gamma $ is simply connected and disjoint from the set $\cup_{j=1}^{n}\rho_{j}(U).$ Theorem \ref{RMT} (the Riemann
  Monodromy theorem)
   gives that each of $\{\rho_{j}(z)\}_{j=1}^{n}$ has a uniquely analytic continuation on  $E\backslash\Gamma $  and hence
   we can view  each of $\{\rho_{j}(z)\}_{j=1}^{n}$ as an analytic function on  $E\backslash\Gamma $ satisfying
  $$\phi (\rho_{j}(z))=\phi (z),$$
  for $z\in E\backslash\Gamma $ and $j=1, \cdots , n.$  We retain
  the same
  labels
   $\{\rho_{j}(z)\}_{j=1}^{n}$ as $\{\rho_{1}(z), \rho_{2}(z), \cdots , \rho_{n}(z)\}$ at every point in $E\backslash\Gamma $.
   The orientation is denoted by $\{\{\rho_{j}(z)\}_{j=1}^{n}, U\}$. Then the composition $\rho_{k}\circ\rho_{l}(z)$
   makes sense on $U$.  As we pointed out before,    $\{\rho_{j}\}_{j=1}^{n}$ has a group-like property under
    composition
   on $U$, that is,
  \begin{equation}\label{kipi}
  \rho_{k}\circ\rho_{i}(z)=\rho_{\pi_{i}(k)}(z),
  \end{equation}
  for some $\pi_{i}(k)\in \{1,2, \cdots , n\}.$  Thus
  \begin{eqnarray}
 \label{piikc} \{\pi_{i}(k)\}_{k=1}^{n}&=&\{1,2, \cdots , n\}.
  \end{eqnarray}
  For each $\rho \in \{\rho_{i}\}_{i=1}^n$, there is a unique
  $\hat{\rho}$ in $\{\rho_{i}\}_{i=1}^n$ such that
  $$\hat{\rho}(\rho (z))=\rho_{1}(z)$$
  on $U$. Then the mapping $\rho\rightarrow \hat{\rho}$ is a
  bijection from the finite set $\{\rho_{i}\}_{i=1}^n$ to itself.
  Thus $\rho_{k}=\hat{\rho}_{k^{-}}$ for some $k^{-}$.
  So for each fixed $k$, there is a unique number $k^{-}$ in $\{1,2, \cdots , n\}$ such that
  $$\rho_{k} (z)=\hat{\rho}_{k^{-}}$$
  and hence
  $$\rho_{k}\circ\rho_{k^{-}}(z)=\rho_{1}(z);$$
  and for each fixed $i$,
  $$\pi_{i}=\left(\begin{array}{cccc}
              1& 2 & \cdots & n\\
              \pi_{i}(1) &\pi_{i}(2) &\cdots &\pi_{i}(n)
              \end{array}\right) $$
              is in   the permutation group $P_{n}.$ Now we define
    a mapping $\Phi$ from  the set  $\{\rho_{j}\}$ of local inverses to  the permutation group $P_{n}$ as
  $$\Phi (i):=\pi_{i}.$$

  Thus for each open set $  \Delta \subset U$,
  \begin{equation}\label{compik}\rho_{k}(\rho_{i}(  \Delta ))=\rho_{\pi_{i}(k)}(  \Delta ).
  \end{equation}

  Note that $\rho_{k}^{-1}(z)$ is also admissible in $E$ and so it is analytic locally in $E\backslash\Gamma$.
   But on $U,$  $\rho_{k}^{-1}(z)=\rho_{k_{-}}(z)$ for some $k_{-}.$
  Thus  \begin{equation}\label{rhoinv}
  \rho_{k}^{-1}(z)=\rho_{k_{-}}(z).
  \end{equation} for $z$ in $E\backslash\Gamma$. So for each subset $  \Delta$ of $E\backslash\Gamma$,
  \begin{equation}\label{rhoinvset}
  \rho_{k}^{-1}(  \Delta )=\rho_{k_{-}}(  \Delta).
  \end{equation}

The proof of Theorem \ref{localrepu} gives a global  form in terms of the orientation defined as above:

\begin{thm}\label{localrep} Let $\phi$ be a finite Blaschke product. Let $U$ be a small invertible  open set of $E$.
Let $\{\rho_{i}(z)\}_{j=1}^{n}$  a complete collection of local inverses  and $E\backslash\Gamma$ with the orientation
$\{\{\rho_{j}(z)\}_{j=1}^{n}, U\}$. Then for  each $T$ in $\{M_{\phi}\}^{\prime}$, there are analytic functions $\{s_{i}(\alpha )\}_{i=1}^{n}$
  on $E\backslash\Gamma$ such that for each $h$ in the Bergman space $L^{2}_{a},$
  \begin{eqnarray*}
  Th(\alpha )&= \sum_{i=1}^{n}s_{i}(\alpha )h(\rho_{i}(\alpha )),\\
  T^{*}k_{\alpha}&=\sum_{i=1}^{n}\overline{s_{i}(\alpha )}k_{\rho_{i}(\alpha )}
  \end{eqnarray*}
 for each $\alpha$ in $E\backslash\Gamma$.
  \end{thm}

Recall that an operator $S$ on a Hilbert space $H$ is subnormal if
there are a Hilbert space $K$ containing $H$ and a normal operator
$N$ on $K$ such that $H$ is an invariant subspace of $N$ and
$S=N|_{H}.$ If $K$ has no proper subspace that contains $H$ and
reduces $N$, we say that $N$ is a minimal normal extension of $S$.
Note that $M_{\phi}$ is a subnormal operator. The main idea is to
use the property that the spectrum of the minimal normal extension
of $M_{\phi}$ contains each invertible set $V$ and the minimal
normal extension of $M_{\phi}$ can be defined using the functional
calculus of the minimal normal extension of the Bergman shift
$M_{z}$ under $\phi$.

We need a few of results on subnormal operators in \cite{Con}:

For $\phi$ in $L^{\infty}$, let $\tilde{M}_{\phi}$ denote
the multiplication operator by $\phi$ on $L^{2}(\D, dA)$ given by
$$\tilde{M}_{\phi}g=\phi g$$
for each $g$ in $L^{2}(\D, dA).$

 \begin{itemize}

\item{\sl
  The minimal normal extension of the Bergman shift $M_{z}$ is the operator $\tilde{M}_{z}$ on $L^{2}(\D, dA)$.}

  \item{\sl For each finite Blaschke product $\phi$, the minimal norm extension of $M_{\phi}$ is the operator $\tilde{M}_{\phi}$ on $L^{2}(\D, dA)$.}

  \item{\sl } Suppose $S_{1}$ and $S_{2}$ are two subnormal operators on a Hilbert space $H$
  and $N_{1}$ and $N_{2}$ are the minimal normal extension of $S_{1}$ and $S_{2}$ on $K$, respectively. If $S_{1}$ and $S_{2}$ are unitarily equivalent,
  i.e., $W^{*}S_{1}W=S_{2}$ for some unitary operator $W$ on $H$, then there is a unitary operator $\tilde{W}$ on $K$ such that
  $$W=\tilde{W}|_{H},$$
  $$\tilde{W}^{*}N_{1}\tilde{W}=N_{2}.$$

\end{itemize}

Thus for each unitary operator $W$ in the commutant of $ M_{\phi}$, there is a unitary operator
$\tilde{W}$ on $L^{2}(\D, dA)$ such that
$$W=\tilde{W}|_{L^{2}_{a}},$$
$$\tilde{W}^{*}\tilde{M}_{\phi}\tilde{W}=\tilde{M}_{\phi}.$$

The following theorem was obtained in \cite{Sun}. For completeness,
we will give the detailed proof of the theorem. A different proof
was also given in \cite{GuH}.

\begin{thm}\label{repuni} Let $\phi$ be a finite Blaschke product.  Let $U$ be a small invertible  open set of $E$. Let $\{\rho_{i}(z)\}_{j=1}^{n}$ be  a complete collection of local inverses  and $E\backslash\Gamma$ with the orientation $\{\{\rho_{j}(z)\}_{j=1}^{n}, U\}$.  Then for   each unitary operator $W$ in $\{M_{\phi}\}^{\prime}$, there is a unit vector  $\{r_{i}\}_{i=1}^{n}$ in $C^{n}$
    such that for each $h$ in the Bergman space $L^{2}_{a},$
\begin{eqnarray*}
  Wh(\alpha )&=& \sum_{i=1}^{n}r_{i}\rho_{i}^{\prime}(\alpha )h(\rho_{i}(\alpha )),\\
  W^{*}k_{\alpha}&=& \sum_{i=1}^{n}\overline{r_{i}}\>\overline{\rho_{i}^{\prime}(\alpha )}k_{\rho_{i}(\alpha )}
  \end{eqnarray*}
for each $\alpha$ in $E\backslash\Gamma.$
\end{thm}

\proof Let $W$ be a unitary operator in the commutant $\{M_{\phi}\}^{\prime}$. By the remark as above,
there is a unitary operator $\tilde{W}$ on $L^{2}(\D, dA)$ such that
$$W=\tilde{W}|_{L^{2}_{a}},$$
$$\tilde{W}^{*}\tilde{M}_{\phi}\tilde{W}=\tilde{M}_{\phi}.$$
Since $\tilde{M}_{\phi}$ is a normal operator on $L^{2}(\D, dA),$ and $\sigma (\tilde{M}_{\phi})=\overline{\D},$
 by the spectral theorem for normal operators \cite{Con}, $\tilde{M}_{\phi}$ has the following spectral decomposition:
$$\tilde{M}_{\phi}=\int_{\overline{\D}}\lambda dE(\lambda ).$$
Note that $\tilde{M}_{\phi}$ equals $\phi (\tilde{M}_{z})$. Thus $\tilde{M}_{\phi}$ commutes with $\tilde{M}_{z}.$ So for each Lebesgue measurable subset $  \Delta$ of $\overline{\D}$, the spectral measure
$$E(  \Delta )=\int_{  \Delta}dE(\lambda )$$
is a projection commuting  with $\tilde{M}_{z}$. Since the commutant
$$\{\tilde{M}_{z}\}^{\prime}=\{\tilde{M}_{f}:f\in L^{\infty}(\D, dA)\}, $$
  we obtain
$$E(  \Delta )=\tilde{M}_{\chi_{\phi^{-1}(  \Delta )}}.$$
Here $\chi_{\phi^{-1}(  \Delta )}$ is the characteristic function of the set $\phi^{-1}(  \Delta ):$
$$\chi_{\phi^{-1}(  \Delta )}(z)=\left\{\begin{array}{cc}
              1&  ~z\in \phi^{-1}(  \Delta )\\
              0 &   ~z\notin\phi^{-1}(  \Delta ) .
              \end{array}\right.
              $$
Since   $\tilde{W}$ commutes with $\tilde{M}_{\phi}$, we have that $E(  \Delta )$ commutes with $\tilde{W}$ to get
$$\tilde{W}^{*}\tilde{M}_{\chi_{\phi^{-1}(  \Delta )}}\tilde{W}=\tilde{M}_{\chi_{\phi^{-1}(  \Delta )}}.$$
Letting $Supp (g)$ denote the essential support
of a function $g$ in $L^{2}(\D, dA),$ the above equality gives
\begin{eqnarray*}
Supp[\tilde{W}^{*}({\chi_{\phi^{-1}(  \Delta )}}g)]&=& Supp [\tilde{W}^{*}\tilde{M}_{\chi_{\phi^{-1}(  \Delta )}}(g)]\\
                                                &=& Supp [\tilde{W}^{*}\tilde{M}_{\chi_{\phi^{-1}(  \Delta )}}\tilde{W} \tilde{W}^{*}(g)]\\
                                                &=&  Supp [\tilde{M}_{\chi_{\phi^{-1}(  \Delta )}}  \tilde{W}^{*}(g)]\\
                                                &=&  Supp [ {\chi_{\phi^{-1}(  \Delta )}}  \tilde{W}^{*}(g)]\\
                                                &\subset & {{\phi^{-1}(  \Delta )}},
                                                \end{eqnarray*}
for each Borel set $    \Delta $ and each $g$ in $L^{2}(\D, dA).$

 For each $z$ in $U$, let $V$ be an open neighborhood of $\alpha$ that is contained in $U$. Since $V$ is a subset of $E$, for each $\alpha$ in $V$, $\phi^{\prime}(\alpha )\neq 0$. Note that on $V$
$$\phi (\rho_{i}(\alpha ))=\phi (\alpha ).$$
Taking the derivative of both sides of the above equality, by the
chain rule we have
$$\phi^{\prime}(\rho_{i}(\alpha ))\rho_{i}^{\prime}(\alpha )=\phi^{\prime}(\alpha )$$
to get that $\rho^{\prime}(\alpha )\neq 0$. This gives that $\rho_{i}$ is locally injective. Since $V$ is a subset of $U$ and $\{\overline{\rho_{i}(U)}\}_{i=1}^{n}$ are disjoint,
$\{\overline{\rho_{i}(V )}\}_{i=1}^{n}$ are disjoint. Thus
$$\phi^{-1}\circ\phi (V )=\cup_{i=1}^{n} \rho_{i}(V)$$
is the union of $n$ strictly separated and polynomially convex open sets where  $\rho_{1}(V)=V$.
Since $\{\overline{\rho_{i}(V)}\}_{i=1}^{n}$ are disjoint,
 $f$ is analytic on the closure of  $\phi^{-1}\circ\phi (V )$ if
$$f(z)=\left\{\begin{array}{ll}
            1 &  z\in \overline{V}\\
            0 & z \in \overline{\phi^{-1}\circ\phi (V )}\backslash\overline{V}.
\end{array} \right.
$$
 By the Runge theorem
\cite{Con}, there exists a sequence of polynomials $\{p_{k}(z)\}$ of $z$  such that $p_{k}(z)$
 uniformly converges to $f(z)$ on the closure of $\phi^{-1}\circ\phi (V)$. Therefore,
\begin{equation}\label{lim1}
\lim_{k\rightarrow \infty}\| \chi_{\phi^{-1}\circ\phi
(V)}p_{k}-\chi_{\phi^{-1}\circ\phi (V)}f\|_{2}=0.
\end{equation}
This gives that for each $g$ in $L^{2}(\D,dA),$
 \begin{eqnarray}
\ip{\tilde{W}^{*}(\chi_{\phi^{-1}\circ\phi (V)}p_{k}), g}
 &\rightarrow & \ip{\tilde{W}^{*}(\chi_{\phi^{-1}\circ\phi (V)}f),  g }\nonumber\\
 \label{Wlimr}&=&\ip{\tilde{W}^{*}( \chi_{V}), g}.
\end{eqnarray}
The last equality follows from the fact that
$$  V\subset \phi^{-1}\circ\phi (V).$$

On the other hand, since $W$ is unitary and commutes with
$M_{\phi}$, $W^{*}$ also commutes with $M_{\phi}$. Thus Theorem
\ref{localrep} gives that there are   functions $t_{i}(\alpha )$ and
$s_{i}(\alpha )$ analytic on $E\backslash\Gamma$  such that
\begin{eqnarray}
Wh(\alpha )&=&\sum_{i=1}^{n}t_{i}(\alpha )h(\rho_{i}(\alpha ))\nonumber,\\
\label{Wstarr}W^{*}h(\alpha ) &=&\sum_{i=1}^{n} {s_{i}(\alpha )}h(\rho_{i}(\alpha ))
\end{eqnarray}
for $\alpha$ in $E\backslash\Gamma$ and $h$ in $L^{2}_{a}.$  It suffices to show that there is a unit vector $(r_{1}, r_{2}, \cdots , r_{n})$ in $C^{n}$ such that
$$s_{i}(\alpha )=r_{i}\rho_{i}^{\prime}(\alpha )$$
on $E\backslash\Gamma$.

 A simple calculation  gives
\begin{eqnarray*}
\ip{\tilde{W}^{*}(\chi_{\phi^{-1}\circ\phi (V)}p_{k}), g}
 &=& \ip{\tilde{W}^{*} \tilde{M}_{\chi_{\phi^{-1}\circ\phi (V)}}\tilde{W}\tilde{W}^{*}(p_{k}),  g }\\
 &=&\ip{  \tilde{M}_{\chi_{\phi^{-1}\circ\phi (V)}} \tilde{W}^{*}(p_{k}),  g }\\
 &=&\int_{ {\phi^{-1}\circ\phi (V)}}\tilde{W}^{*}(p_{k})(\alpha )\overline{g(\alpha )}dA(\alpha )\\
 &=&\int_{D}\chi_{\phi^{-1}\circ\phi (V)}\sum_{i=1}^{n} s_{i}(\alpha )p_{k}(\rho_{i}(\alpha ))\overline{g(\alpha )}dA(\alpha ).
\end{eqnarray*}
The last equality follows from (\ref{Wstarr}).
Since the closure of $\phi^{-1}\circ\phi (V)$ is a compact subset of $E\backslash\Gamma $ and $\{s_{i}(\alpha )\}_{i=1}^{n}$ are analytic in $E\backslash\Gamma $, there is a constant $M>0$ such that
$$\sup_{1\leq i\leq n}\sup_{\alpha \in \phi^{-1}\circ\phi (V)}|s_{i}(\alpha )|\leq M.$$
Noting that $\rho_{i}^{\prime}(\alpha )\neq 0$ for $\alpha$ in
$E\backslash\Gamma$ and $p_{k}(z)$ uniformly converges to $f(z)$ on
the closure of $\phi^{-1}\circ\phi (V)$,  by (\ref{rhoinvset}), we
see that
$$\lim_{k\rightarrow \infty}\|\chi_{\phi^{-1}\circ\phi (V)}p_{k}\circ\rho_{i}-\chi_{\phi^{-1}\circ\phi (V)}f\circ\rho_{i}\|_{2}=0.$$
Therefore, we have
$$\lim_{k\rightarrow \infty}\|\chi_{\phi^{-1}\circ\phi (V)}\sum_{i=1}^{n} s_{i}p_{k}\circ\rho_{i}   -\chi_{\phi^{-1}\circ\phi (V)}\sum_{i=1}^{n} s_{i} f\circ\rho_{i}  \|_{2}=0,$$
to obtain
\begin{eqnarray}\label{Wliml}
\ip{\tilde{W}^{*}(\chi_{\phi^{-1}\circ\phi (V)}p_{k}), g} &\rightarrow &\ip{\chi_{\phi^{-1}\circ\phi (V)}
\sum_{i=1}^{n} s_{i}(\alpha )f(\rho_{i}(\alpha )), g}.
\end{eqnarray}
Combining (\ref{Wlimr}) and (\ref{Wliml}) gives that
\begin{equation}
\tilde{W}^{*}( \chi_{V})= \chi_{\phi^{-1}\circ\phi
(V)}\sum_{i=1}^{n} s_{i}(\alpha )f(\rho_{i}(\alpha )).
\end{equation}
Noting
  \begin{eqnarray*}
Supp [\sum_{i=1}^{n} s_{i}(\alpha )f(\rho_{i}(\alpha ))]&\subset& \cup_{i=1}^{n} Supp[f\circ\rho_{i}]\\
                                                        &\subset &\cup_{i=1}^{n} Supp[\chi_{\rho_{i}^{-1}(V)}]\\
                                                        &\subset & \cup_{i=1}^{n} \rho_{i^{-}}(V)\hspace{2cm}\text{by}~~~~ (\ref{rhoinvset})\\
                                                        &=&       \phi^{-1}\circ\phi (V),
                                                        \end{eqnarray*}
 we obtain
$$\tilde{W}^{*}( \chi_{V})= \sum_{i=1}^{n} s_{i}(\alpha )f(\rho_{i}(\alpha )).$$
For an area-measurable set $V$ in the unit disk, we use $|V|$ to
denote the area measure of $V$. Since $W$ is a unitary operator on
$L^{2}(\D, dA)$, we have
\begin{eqnarray*}
|V|^{2} &=& \|\chi_{V}\|_{2}^{2}\\
        &=& \|W^{*}(\chi_{V})\|_{2}^{2}\\
        &=& \int_{D}|\sum_{i=1}^{n} s_{i}(\alpha )f(\rho_{i}(\alpha ))|^{2}dA(\alpha )\\
        &=& \int_{D}\sum_{i=1}^{n}\sum_{j=1}^{n}s_{i}(\alpha )\overline{s_{j}(\alpha )}f(\rho_{i}(\alpha ))\overline{f(\rho_{j}(\alpha ))}dA(\alpha )\\
        &=& \int_{D}\sum_{i=1}^{n}|s_{i}(\alpha )|^{2}|f(\rho_{i}(\alpha ))|^{2}dA(\alpha )\\
        &=& \sum_{i=1}^{n}\int_{\rho_{i}^{-1}(V)} |s_{i}(\alpha )|^{2}|f(\rho_{i}(\alpha ))|^{2}dA(\alpha )\\
        &=& \sum_{i=1}^{n}\int_{V}|f(\alpha )|^{2}\Bigl|\frac{s_{i}\circ\rho_{i}^{-1}(\alpha )}{\rho_{i}^{\prime}\circ\rho_{i}^{-1}(\alpha )}\Bigr|^{2}dA(\alpha )\\
        &=& \int_{V}\sum_{i=1}^{n}\Bigl|\frac{s_{i}\circ\rho_{i}^{-1}(\alpha )}{\rho_{i}^{\prime}\circ\rho_{i}^{-1}(\alpha )}\Bigr|^{2}dA(\alpha ).
        \end{eqnarray*}
        The fifth equality follows from
        \begin{eqnarray*}
        Supp[f\circ \rho_{i}]\cap Supp[f\circ \rho_{j}]&=&\rho_{i}^{-1}(V)\cap \rho_{i}^{-1}(V)\\
                                                       &=&\rho_{i_{-}}(V)\cap \rho_{j_{-}}(V)\hspace{2cm}\text{by}~~~(\ref{rhoinvset})\\
                                                       &=& \emptyset.
                                                       \end{eqnarray*}
         Here $\emptyset$ denotes the empty set.
         The sixth equality follows from that $Supp[f\circ\rho_{i}]=\rho_{i}^{-1}(V)$  and
         the seventh equality follows from that
         the change of variable, $\beta =\rho_{i}(\alpha ),$ gives
         $$\int_{\rho_{i}^{-1}(V)} |s_{i}(\alpha )|^{2}|f(\rho_{i}(\alpha ))|^{2}dA(\alpha )=
          \int_{  V } |f(\beta )|^{2} \Bigl|\frac{s_{i}\circ\rho_{i}^{-1}(\beta )}{\rho_{i}^{\prime}\circ\rho_{i}^{-1}(\beta )}\Bigr|^{2}dA(\beta ).$$
          Therefore we obtain  for each $z\in U$ and for any open neighborhood $V\subset U$ of
          $z$, that
                                                         $$1=\frac{1}{{|V|^{2}}}\int_{V}\sum_{i=1}^{n}\Bigl|\frac{s_{i}\circ\rho_{i}^{-1}(\alpha )}{\rho_{i}^{\prime}\circ\rho_{i}^{-1}(\alpha )}\Bigr|^{2}dA(\alpha ).$$
        Noting that $\sum_{i=1}^{n}\Bigl|
        \frac{s_{i}\circ\rho_{i}^{-1}(\alpha )}{\rho_{i}^{\prime}\circ\rho_{i}^{-1}(\alpha )}\Bigr|^{2}$
        is continuous on $U$ and  letting $V$ shrink to $z$, we have
                                                        $$1=\sum_{i=1}^{n}\Bigl|\frac{s_{i}\circ\rho_{i}^{-1}(z )}{\rho_{i}^{\prime}\circ\rho_{i}^{-1}(z )}\Bigr|^{2}.$$
                                                        Applying the Laplace operator $\frac{\partial^{2}}{\partial z\partial \overline{z}}$ to
                                                        both sides of the above equality gives
           $$0=\sum_{i=1}^{n}\Bigl|\Bigl[\frac{s_{i}\circ\rho_{i}^{-1}(z )}{\rho_{i}^{\prime}\circ\rho_{i}^{-1}(z )}\Bigr]^{\prime}\Bigr|^{2}.$$
           Thus
           $$\Bigl[\frac{s_{i}\circ\rho_{i}^{-1}(z )}{\rho_{i}^{\prime}\circ\rho_{i}^{-1}(z )}\Bigr]^{\prime}=0$$
           on $U$. So there are constants $r_{i}$ such that
           \begin{eqnarray*}
           |r_{1}|^{2}+\cdots +|r_{n}|^{2}&=&1\\
           \frac{s_{i}\circ\rho_{i}^{-1}(z )}{\rho_{i}^{\prime}\circ\rho_{i}^{-1}(z )}&=&r_{i}
           \end{eqnarray*}
           on $U$ for each $i$, and so
           $$s_{i}(z )=r_{i}\rho_{i}^{\prime}(z ),$$
           on $\rho_{i}^{-1}(U)=\rho_{i_{-}}(U)$ and hence on $E\backslash\Gamma$ because that both $s_{i}(z )$ and $\rho_{i}^{\prime}(z )$ are analytic on $E\backslash\Gamma$. This completes the proof.

\section{\label{DecompBergman}Decomposition of the Bergman space}

Let $\phi$ be a Blaschke product of order $n$. As pointed out in the
introduction, by complex geometry, one can easily see  that the
Bergman space can be decomposed as a direct sum of at most $n$
nontrivial   reducing subspaces of $M_{\phi}$. In this section we
will refine this result and show that the Bergman space can be
decomposed as a direct sum of at most $q$ nontrivial  reducing
subspaces of $M_{\phi}$, where $q$ is the number of connected
components of the Riemann surface $S_{\phi}$.

Let $U$ be a small invertible  open subset of $E$. Let
$\{\rho_{i}(z)\}_{j=1}^{n}$  be a complete collection of local
inverses. Theorem \ref{repuni} gives that for each $z$ in $U$,
  \begin{eqnarray*}
    W^{*}k_{\alpha}&=& \sum_{i=1}^{n}\overline{r_{i}}\>\overline{\rho_{i}^{\prime}(\alpha )}k_{\rho_{i}(\alpha )}.
  \end{eqnarray*}
  For our convenience, we use $r_{\rho_{i}}$ to denote $r_{i}.$  Then the    above equality becomes
  \begin{eqnarray*}
    W^{*}k_{\alpha}&=& \sum_{\rho\in \{\rho_{i}(z)\}_{j=1}^{n}}\overline{r_{\rho}}\>\overline{\rho^{\prime}(\alpha )}k_{\rho(\alpha )}.
  \end{eqnarray*}

  \begin{lem}\label{decomp} Let $U$ be a small invertible  open subset of
  $E$ and
  $\{\rho_{i}(z)\}_{j=1}^{n}$ be
  a complete collection of local inverses. Fix a point $z_{0}$ in $U$. Let $G_{i}$ be the set of those functions in $\{\rho_{i}(z)\}_{i=1}^{n}$ which are analytic extensions of $\rho_{i}$ along some loop  containing $z_{0}$ in $E$. Write
$$ \{\rho_{i}(z)\}_{i=1}^{n}=\cup_{k=1}^{q}G_{i_{k}}.$$ Then for each $\alpha$ in $U$,
  \begin{eqnarray*}
    W^{*}k_{\alpha}&=& \sum_{k=1}^{q}\hat{r}_{k}\sum_{\rho\in G_{i_{k}}}\overline{\rho^{\prime}(\alpha )}k_{\rho (\alpha )},
  \end{eqnarray*}
  where $\hat{r}_{k}=r_{\rho}$ for some $\rho$ in $G_{i_{k}}.$
  \end{lem}

  \proof
  For each   $\alpha$ in $U$,  by Theorem \ref{repuni},  we have
  \begin{eqnarray}
   \label{Wstarproof} W^{*}k_{\alpha}&=& \sum_{k=1}^{q} \sum_{\rho\in G_{i_{k}}}\overline{r_{\rho}}\> \overline{\rho^{\prime}(\alpha )}k_{\rho (\alpha )}.
  \end{eqnarray}
  It suffices to show that for two $\rho_{0}$ and $\hat{\rho}_{0}$ in the same
  $G_{i_{k}}$, $r_{\rho_{0}}=r_{\hat{\rho}_{0}}$.
  Note that the conjugates of both sides of (\ref{Wstarproof}) are locally analytic
  functions of $\alpha$ in $E$. In fact the conjugate of
  the right hand side of (\ref{Wstarproof}) is   an analytic function of $\alpha$
  in the unit disk $\D$.
 Since each of  $\{\rho_{i}(z)\}_{j=1}^{n}$  admits  unrestricted continuation in $E$, the conjugate
 of the left hand of    (\ref{Wstarproof})    admits  unrestricted analytic continuation in $E$. Let $z_{0}$ be the fixed point in $U$ and $\gamma$
 be   a loop in $E$ containing $z_{0} $.

 Suppose that each $\rho$ of the set $\{\rho_{i}(z)\}_{i=1}^{n}$ is
extended analytically to $\hat{\rho}$ from $z_{0}$ to $z_{0}$ along
the loop $\gamma$.  In the neighborhood $U$ of $z_{0}$, we have
  \begin{eqnarray} \label{Wstarproof1}
  W^{*}k_{\alpha}&=& \sum_{k=1}^{q} \sum_{\rho\in G_{i_{k}}}\overline{r_{\rho}}\> \overline{\hat{\rho}^{\prime}(\alpha )}k_{\hat{\rho} (\alpha )}.
 \end{eqnarray}
 Thus (\ref{Wstarproof}) and (\ref{Wstarproof1}) give
   \begin{eqnarray*}
  \sum_{k=1}^{q} \sum_{\rho\in G_{i_{k}}}\overline{r_{\rho}}\> \overline{\rho^{\prime}(\alpha )}k_{\rho (\alpha )}&=& \sum_{k=1}^{q} \sum_{\rho\in G_{i_{k}}}\overline{r_{\rho}}\> \overline{\hat{\rho}^{\prime}(\alpha )}k_{\hat{\rho} (\alpha )}.
 \end{eqnarray*}
 Since  the reproducing kernels are linearly independent, the above equality gives
 $$\overline{ r_{\hat{\rho}_{0}}}\> \overline{\hat{\rho}^{\prime}_{0}(\alpha )}= \overline{r_{\rho_{0}}}\>
 \overline{\hat{\rho}^{\prime}_{0}(\alpha )}.$$  And hence $r_{\hat{\rho}_{0}}= r_{\rho_{0}}.$ This completes the proof.

 \begin{thm}\label{Bergdecomp} Let $\phi$ be a finite  Blaschke product. Then the Bergman space can be decomposed as a direct
 sum of at most $q$ nontrivial minimal reducing subspaces of $M_{\phi} $ where $q$ is the number of connected components of the Riemann surface $S_{\phi}$.
  \end{thm}

  \proof Suppose that the Bergman space is the direct sum of  $p$ nontrivial minimal reducing subspaces $\{{\mathcal M}_{j}\}_{j=1}^{p}$ of $M_{\phi}$. That is,
  $$L^{2}_{a}=\oplus_{j=1}^{p}{\mathcal M}_{j} .$$
  Let $P_{j}$ denote the orthogonal projection from $L^{2}_{a}$ onto ${\mathcal M}_{j}$. Thus $P_{j}$ commutes with
  both $M_{\phi}$ and $M_{\phi}^{\prime}$. For $(\theta_{1}, \cdots , \theta_{p})$  in $[0,2\pi]^{n}$, let
  $$W(\theta_{1}, \cdots , \theta_{p})=\sum_{j=1}^{p}e^{i\theta_{j}}P_{j}.$$
  Then $\{W(\theta_{1}, \cdots , \theta_{p})\}$ is a family of unitary operators in the commutant
  $\{M_{\phi}\}^{\prime}.$ Using this family of unitary operators, we can recover $P_{j}$ as follows:
  $$P_{j}h=\int_{[0,2\pi]^{p}}[W(\theta_{1}, \cdots , \theta_{p})]^{*}hd\mu_{j}(\theta_{1}, \cdots , \theta_{p}),$$
  for $h$ in $L^{2}_{a},$ where $d\mu_{j}(\theta_{1}, \cdots , \theta_{p})$ is the measure $\frac{e^{i\theta_{j}}}{(2\pi )^{p}}d\theta_{1}\cdots\theta_{p}$
  on $[0,2\pi]^{p}.$

  Let $U$ be a small invertible  open subset of $E$ and $\{\rho_{i}(z)\}_{j=1}^{n}$ be  a complete collection of
  local inverses. Fix a point $z_{0}$ in $U$. Let $G_{i}$ be the set of those local inverses
  in the set
  $\{\rho_{i}(z)\}_{i=1}^{n}$
  which are analytic extensions of $\rho_{i}$ along some loop   in $E$ containing $z_{0}$. Write
$ \{\rho_{i}(z)\}_{i=1}^{n}=\cup_{k=1}^{q}G_{i_{k}}.$ Theorem
\ref{Riem} shows that $q$ is the number   of connected components of
the Riemann surface $S_{\phi}$.   To finish the proof we need only
to show that
$$p\leq q.$$
To do this, by Lemma \ref{decomp}, we have that  there are constants
  $\{\hat{r}_{k}(\theta_{1}, \cdots , \theta_{p})\}$ such that
  $$[W(\theta_{1}, \cdots , \theta_{p})]^{*}k_{\alpha}= \sum_{k=1}^{q}\hat{r}_{k}(\theta_{1}, \cdots , \theta_{p})\sum_{\rho\in G_{i_{k}}}\overline{\rho^{\prime}(\alpha )}k_{\rho (\alpha )}$$
  for $\alpha$ in $U$. Thus
  \begin{eqnarray}
  P_{j}k_{\alpha}&=& \int_{[0,2\pi]^{p}}[W(\theta_{1}, \cdots , \theta_{p})]^{*}k_{\alpha}d\mu_{j}(\theta_{1}, \cdots , \theta_{p}) \nonumber\\
                 &=& \int_{[0,2\pi]^{p}}\sum_{k=1}^{q}\hat{r}_{k}(\theta_{1}, \cdots , \theta_{p})\sum_{\rho\in G_{i_{k}}}\overline{\rho^{\prime}(\alpha )}k_{\rho (\alpha )}d\mu_{j}(\theta_{1}, \cdots , \theta_{p})\nonumber\\
               \label{PjWj}  &=& \sum_{k=1}^{q}\tilde{r}_{k} \sum_{\rho\in G_{i_{k}}}\overline{\rho^{\prime}(\alpha )}k_{\rho (\alpha )}
                 \end{eqnarray}
                 where
                 $$\tilde{r}_{k}=\int_{[0,2\pi]^{p}}\hat{r}_{k}(\theta_{1}, \cdots , \theta_{p})d\mu_{j}(\theta_{1}, \cdots , \theta_{p}).$$
                 For each $j$, let ${\mathcal Z}_{j}$ denote the zero set of the functions in ${\mathcal
                 M}_{j}$, that is
                 $$\{z\in \D: f(z)=0~~ \text{for each }~~ f\in {\mathcal M}_{j}\}.$$
                 Then ${\mathcal Z}_{j}$ is a countable subset of the unit disk $\mathbb D$. Hence $U\backslash[\cup_{j=1}^{p}{\mathcal Z}_{j}]$
                 is not empty. For each $\alpha$ in $U\backslash[\cup_{j=1}^{p}{\mathcal Z}_{j}]$,
                 there is a function $f_{j}$ in ${\mathcal M}_{j}$
                 such that
                 $$f_{j}(\alpha )\neq 0.$$
                 Further, we have
                 \begin{eqnarray*}
                 f_{j}(\alpha )&=&\ip{f_{j}, k_{\alpha}}\\
                               &=& \ip{P_{j}f_{j}, k_{\alpha}}\\
                               &=& \ip{f_{j}, P_{j}k_{\alpha}}.\\
                 \end{eqnarray*}
                 Thus
                 \begin{equation}\label{mjspace}
                 \|P_{j}k_{\alpha}\|\neq 0.
                 \end{equation}
                 For a fixed $\alpha$ in $U\backslash[\cup_{j=1}^{p}{\mathcal Z}_{j}]$, let $\Theta$ be the subspace of $L^{2}_{a}$  spanned by $\{P_{1}k_{\alpha}, \cdots P_{p}k_{\alpha}\}.$ Thus (\ref{PjWj})  gives that
                 $\Theta$ is contained in the subspace of $L^{2}_{a}$ spanned by $q$ functions
                 $$\{\sum_{\rho\in G_{i_{k}}}\overline{\rho^{\prime}(\alpha )}k_{\rho (\alpha )}\}_{k=1}^{q}.$$  So the dimension of $\Theta$ is less than
                 or equal to $q$.
                 On the other hand, for distinct $j$ and $l$, $P_{l}P_{j}=0$. This gives
                 \begin{eqnarray*}
                 \ip{P_{j}k_{\alpha}, P_{l}k_{\alpha}}&=& \ip{P_{l}P_{j}k_{\alpha}, k_{\alpha}}\\
                                                      &=& 0.
                                                      \end{eqnarray*}
                  Thus combining the above equality with (\ref{mjspace}) gives that   the dimension of $\Theta$ equals $p$ and so $p\leq q$. We  complete the proof.

\section{matrix representation of unitary operators}

For    a finite Blaschke product $\phi$,  it is pointed out in
Section \ref{repUO} that the orientation of $\{\rho_{i}\}_{i=1}^{n}$
 induces elements $\{\pi_{i}\}_{i=1}^{n}$ in $P_{n}$. So
$\{\pi_{i}\}_{i=1}^{n}$ forms a subgroup  of the permutation group
$P_{n}$. In this section, using the property that
$\{\pi_{i}\}_{i=1}^{n}$ forms a group,  we will obtain the following
unitary matrix representation of a unitary operator in the commutant
$\{M_{\phi}\}^{\prime}$.


\begin{thm}\label{rep} Let $\phi$ be a finite Blaschke product with order $n$, $U$ be a small invertible  open set of
$E$, and  $\{\rho_{j}(z)\}_{j=1}^{n}$  be a complete collection of
local inverses. For $E\backslash\Gamma$ with the orientation
$\{\{\rho_{j}(z)\}_{j=1}^{n}, U\}$. If  $W$ is a unitary operator on
the Bergman space which  commutes with the multiplication operator
$M_{\phi}$, then there is a unit vector $(r_{1}, \cdots, r_{n})\in
{\mathbb C}^{n}$  such that
\begin{equation}\label{matrxrep}W^{*}\left(\begin{array}{c}
   \overline{\rho_{1}^{\prime}(\alpha )}k_{\rho_{1}(\alpha )}\\
\overline{\rho_{2}^{\prime}(\alpha )}k_{\rho_{2}(\alpha )}\\
\overline{\rho_{3}^{\prime}(\alpha )}k_{\rho_{3}(\alpha )}\\
\vdots\\
\overline{\rho_{n}^{\prime}(\alpha )}k_{\rho_{n}(\alpha )}\end{array}\right)=\Gamma_{W}
\left(\begin{array}{c}
   \overline{\rho_{1}^{\prime}(\alpha )}k_{\rho_{1}(\alpha )}\\
\overline{\rho_{2}^{\prime}(\alpha )}k_{\rho_{2}(\alpha )}\\
\overline{\rho_{3}^{\prime}(\alpha )}k_{\rho_{3}(\alpha )}\\
\vdots\\
\overline{\rho_{n}^{\prime}(\alpha )}k_{\rho_{n}(\alpha )}\end{array}\right),
\end{equation}
for $\alpha \in U, $ where the representing matrix $\Gamma_{W}$ is
the following unitary matrix in $U_{n}({\mathbb C})$:
$$\Gamma_{W}=\left(\begin{array}{ccccc}
r_{1}&r_{2}&r_{3}&\cdots&r_{n}\\
r_{\pi_{2}^{-1}(1)}&r_{\pi_{2}^{-1}(2)}&r_{\pi_{2}^{-1}(3)}&\cdots&r_{\pi_{2}^{-1}(n)}\\
r_{\pi_{3}^{-1}(1)}&r_{\pi_{3}^{-1}(2)}&r_{\pi_{3}^{-1}(3)}&\cdots&r_{\pi_{3}^{-1}(n)}\\
\vdots&\vdots&\vdots&\cdots&\vdots\\
r_{\pi_{n}^{-1}(1)}&r_{\pi_{n}^{-1}(2)}&r_{\pi_{n}^{-1}(3)}&\cdots&r_{\pi_{n}^{-1}(n)} \end{array}\right).
$$

\end{thm}

\proof By Theorem \ref{repuni},  there is a unit vector
$(r_{1}, r_{2}, \cdots , r_{n})$ in $C^{n}$ such that
\begin{eqnarray}\label{Wstar1}
W^{*}k_{\alpha}&=& \sum_{i=1}^{n}\overline{r_{i}}\>\overline{\rho_{i}^{\prime}(\alpha )}k_{\rho_{i}(\alpha )}
  \end{eqnarray}
for each $\alpha$ in $E\backslash\Gamma.$ Since
$\cup_{j=1}^{n}\rho_{j}(U)$ is contained in $E\backslash\Gamma$, we
have  for each $l$ and $\alpha$ in $U$, that
\begin{eqnarray}
W^{*}[\overline{\rho_{l}^{\prime}(\alpha )}k_{\rho_{l}(\alpha )}]&=& \overline{\rho_{l}^{\prime}(\alpha )} \sum_{i=1}^{n}\overline{r_{i}}\>\overline{\rho_{i}^{\prime}(\rho_{l}(\alpha ))}k_{\rho_{i}(\rho_{l}(\alpha ))} \nonumber\\
&=&\sum_{i=1}^{n}\overline{r_{i}}\>\overline{\rho_{i}^{\prime}(\rho_{l}(\alpha ))}\>\overline{\rho_{l}^{\prime}(\alpha )}k_{\rho_{i}(\rho_{l}(\alpha ))}\nonumber\\
&=&\sum_{i=1}^{n}\overline{r_{i}}\>\overline{[\rho_{i}\circ\rho_{l}]^{\prime}(\alpha )}k_{\rho_{i}\circ\rho_{l}(\alpha )}\nonumber\\
&=&\sum_{i=1}^{n}\overline{r_{i}}\>\overline{\rho_{\pi_{l}(i)}^{\prime}(\alpha )}k_{\rho_{\pi_{l}(i)}(\alpha )}
\hspace{2cm}\text{by}~~~ (\ref{kipi})\nonumber\\
\label{Wstarl}&=&\sum_{i=1}^{n}\overline{r_{\pi_{l}^{-1}(i)}}\>\overline{\rho_{i}^{\prime}(\alpha )}k_{\rho_{i}(\alpha )}.
  \end{eqnarray}
  This gives the matrix representation (\ref{matrxrep}) of $W$. To finish the proof we need only to show that
  $\Gamma_{W}$ is a unitary matrix.

By  Fuglede Theorem \cite{Con}, \cite{Dou}, $W^{*}$ also commutes
with $M_{\phi}$. For a unitary operator $W$, this follows easily from
the following reason:
$$
WM_{\phi}=M_{\phi}W.
$$
 By multiplying   both sides of the above equality by $W^{*}$  and noting that $$W^{*}W=WW^{*}=I,$$ we have
$$M_{\phi}W^{*}=W^{*}M_{\phi} .$$

Using Theorem \ref{repuni} again,  there is another unit vector
$(s_{1}, s_{2}, \cdots , s_{n})$ in $C^{n}$ such that
\begin{eqnarray}\label{Wl}
Wk_{\alpha}&=& \sum_{i=1}^{n}\overline{s_{i}}\>\overline{\rho_{i}^{\prime}(\alpha )}k_{\rho_{i}(\alpha )}
  \end{eqnarray}
for each $\alpha$ in $E\backslash\Gamma.$  Similar to the argument
estimating (\ref{Wstarl}), we also have
\begin{eqnarray*}
W[\overline{\rho_{l}^{\prime}(\alpha )}k_{\rho_{l}(\alpha )}]&=& \sum_{i=1}^{n}\overline{s_{\pi_{l}^{-1}(i)}}\>\overline{\rho_{i}^{\prime}(\alpha )}k_{\rho_{i}(\alpha )},
  \end{eqnarray*} for $\alpha$ in $U$. Thus for each $l$, we have
\begin{eqnarray*}
\overline{\rho_{l}^{\prime}(\alpha )}k_{\rho_{l}(\alpha )} &=& WW^{*}[\overline{\rho_{l}^{\prime}(\alpha )}k_{\rho_{l}(\alpha )}]\\
           &=& W[W^{*}(\overline{\rho_{l}^{\prime}(\alpha )}k_{\rho_{l}(\alpha )})]\\
           &=&W [\sum_{i=1}^{n}\overline{r_{\pi_{l}^{-1}(i)}}\>\overline{\rho_{i}^{\prime}(\alpha )}k_{\rho_{i}(\alpha )}]\hspace{2cm}\text{by}~~ (\ref{Wstarl})\\
           &=& \sum_{i=1}^{n}\overline{r_{\pi_{l}^{-1}(i)}}\>W[\overline{\rho_{i}^{\prime}(\alpha )}k_{\rho_{i}(\alpha )}]\\
           &=& \sum_{i=1}^{n}\overline{r_{\pi_{l}^{-1}(i)}}\>\sum_{j=1}^{n}\overline{s_{\pi_{i}^{-1}(j)}}\>\overline{\rho_{j}^{\prime}(\alpha )}k_{\rho_{j}(\alpha )}\\
 &=&\sum_{j=1}^{n}\Bigl(\sum_{i=1}^{n} \overline{r_{\pi_{l}^{-1}(i)}}\>\>
         \overline{s_{\pi_{i}^{-1}(j)}}\Bigr)\overline{\rho_{j}^{\prime}(\alpha )}k_{\rho_{j}(\alpha )}.
           \end{eqnarray*}
         Since the reproducing kernels are linearly independent, we have
   \begin{eqnarray}\label{rs1}
   \sum_{i=1}^{n}\overline{r_{\pi_{l}^{-1}(i)}}\>\overline{s_{\pi_{i}^{-1}(l)}}&=&1\\
   \label{rsj}\sum_{i=1}^{n}\overline{r_{\pi_{l}^{-1}(i)}}\>\overline{s_{\pi_{i}^{-1}(j)}}&=&0
    \end{eqnarray}
    for $1\leq j, l\leq n$ and  $j\neq l$. On the other hand,   since  both $(r_{1}, r_{2}, \cdots , r_{n})$  and  $(s_{1}, s_{2}, \cdots , s_{n})$ are unit vectors in $C^{n}$,  by (\ref{piikc})
     we have that for each $k$
   \begin{eqnarray*}
    |r_{\pi_{k}^{-1}(1)}|^{2}+ |r_{\pi_{k}^{-1}(2)}|^{2} + \cdots +|r_{\pi_{k}^{-1}(n)}|^{2} &=&  |r_{1}|^{2}+ |r_{2}|^{2} + \cdots +|r_{n}|^{2}\\
                                                 &=& 1,~~\text{and}\\
     \sum_{k=1}^{n}[|s_{\pi_{1}^{-1}(k)}|^{2}+|s_{\pi_{2}^{-1}(k)}|^{2}+ \cdots +|s_{ \pi_{n}^{-1}(k)}|^{2}]
     &=&\sum_{j=1}^{n}[|s_{\pi_{j}^{-1}(1)}|^{2}+|s_{\pi_{j}^{-1}(2)}|^{2}+ \cdots +|s_{ \pi_{j}^{-1}(n)}|^{2}]\\
     &=& \sum_{j=1}^{n}[|s_1|^{2}+|s_{2}|^{2}+ \cdots +|s_{n}|^{2}]= n.
                              \end{eqnarray*}
                              Let $|s_{\pi_{1}^{-1}(l)}|^{2}+|s_{\pi_{2}^{-1}(l)}|^{2}+ \cdots +|s_{
                              \pi_{n}^{-1}(l)}|^{2}$ be the smallest
                                 of the set
                              $$\{|s_{\pi_{1}^{-1}(k)}|^{2}+|s_{\pi_{2}^{-1}(k)}|^{2}+ \cdots +|s_{
                              \pi_{n}^{-1}(k)}|^{2}\}_{k=1}^{n}.$$
                             Then the pigeonhole principle implies
                             that
                              $$|s_{\pi_{1}^{-1}(l)}|^{2}+|s_{\pi_{2}^{-1}(l)}|^{2}+ \cdots +|s_{
                              \pi_{n}^{-1}(l)}|^{2}\leq 1.$$
                            From the  Cauchy-Schwarz inequality,
                            (\ref{rs1}), one has
                             \begin{eqnarray*}
                             1&=&|\sum_{i=1}^{n}\overline{r_{\pi_{l}^{-1}(i)}}\>\overline{s_{\pi_{i}^{-1}(l)}}|\\
     &\leq & \Bigl(\sum_{i=1}^{n}|r_{\pi_{l}^{-1}(i)}|^{2}\Bigr)^{1/2}\Bigl(\sum_{i=1}^{n}|s_{\pi_{i}^{-1}(l)}|^{2}\Bigr)^{1/2}\leq
     1,
            \end{eqnarray*}
 which gives
     \begin{eqnarray*}(r_{\pi_{l}^{-1}(1)} , r_{\pi_{l}^{-1}(2)} , \cdots , r_{\pi_{l}^{-1}(n)} )&=&
     \lambda_l (\overline{s}_{\pi_{1}^{-1}(l)}, \overline{s}_{\pi_{2}^{-1}(l)}, \cdots , \overline{s}_{ \pi_{n}^{-1}(l)}),
     \end{eqnarray*}
     for some unimodular constant $\lambda_l$. Thus
$$|s_{\pi_{1}^{-1}(l)}|^{2}+|s_{\pi_{2}^{-1}(l)}|^{2}+ \cdots +|s_{
                              \pi_{n}^{-1}(l)}|^{2}=1,$$
                              so we have
                              $$\sum_{k=1, k\neq l}^{n}[|s_{\pi_{1}^{-1}(k)}|^{2}+|s_{\pi_{2}^{-1}(k)}|^{2}+ \cdots +|s_{ \pi_{n}^{-1}(k)}|^{2}]
=n-1.$$ Repeating the above argument and by induction, we will
obtain that for each $k$
\begin{eqnarray*}(r_{\pi_{k}^{-1}(1)} , r_{\pi_{k}^{-1}(2)} , \cdots , r_{\pi_{k}^{-1}(n)} )&=&
     \lambda_k (\overline{s}_{\pi_{1}^{-1}(k)}, \overline{s}_{\pi_{2}^{-1}(k)}, \cdots , \overline{s}_{ \pi_{n}^{-1}(k)}),
     \end{eqnarray*}
     for a unimodular constant $\lambda_k$.

     By  (\ref{rsj}), we have
    \begin{eqnarray*}
    (r_{\pi_{l}^{-1}(1)} , r_{\pi_{l}^{-1}(2)} , \cdots , r_{\pi_{l}^{-1}(n)} )&\perp& (\overline{s}_{\pi_{1}^{-1}(j)}, \overline{s}_{\pi_{2}^{-1}(j)}, \cdots , \overline{s}_{ \pi_{n}^{-1}(j)})
    \end{eqnarray*} for $j\neq l$.
    So \begin{eqnarray*}
    |r_{\pi_{l}^{-1}(1)}|^{2}+ |r_{\pi_{l}^{-1}(2)}|^{2} + \cdots +|r_{\pi_{l}^{-1}(n)}|^{2} &=&1~~\text{for $1\leq l\leq n$ and}\\
    (r_{\pi_{l}^{-1}(1)} , r_{\pi_{l}^{-1}(2)} , \cdots , r_{\pi_{l}^{-1}(n)} )&\perp& (r_{\pi_{j}^{-1}(1)}, r_{\pi_{j}^{-1}(2)}, \cdots , r_{ \pi_{j}^{-1}(n)})
    \end{eqnarray*} for $1\leq j\neq l\leq n$. This gives that $\Gamma_{W}$ is unitary to complete the proof.

\section{   von Neumann algebra ${\mathcal A}_{\phi}$  \label{mainsection}}

In the previous section, we have shown that ${\mathcal A}_{\phi}$ is
finite dimensional. In this section we will show that its dimension
equals the number $q$ of connected components of the Riemann surface
$\phi^{-1}\circ\phi$ over the unit disk. Since Lemma \ref{decomp}
will show directly that the dimension of ${\mathcal A}_{\phi}$ is no
greater than $q$, the main effort in this section is showing that
the dimension is at least $q$.

To that end we are going to construct $q$ linearly independent
elements in
 ${\mathcal A}_{\phi}$. Let $\phi$ be a finite Blaschke product.
Recall that $ {\mathcal C}$ denotes
 the set of the critical points of $\phi$ in $ \D$ and
 $$ {\mathcal F}=\phi^{-1}\circ\phi ( {\mathcal C}),$$
and  $E={ \D}\backslash {\mathcal F}$. Let $z_{0}$ be a point in
 a small invertible  open set $U$ of $E$. Let
$\{\rho_{i}(z)\}_{i=1}^{n}$ be a complete collection of local
inverses.   Let $G_{i}$ be the set of those functions in
$\{\rho_{i}(z)\}_{i=1}^{n}$ which are analytic extensions of
$\rho_{i}$ along some loop  in $E$ containing $z_{0}$. Then Theorem
\ref{Riem}  gives
$$
 \{\rho_{i}(z)\}_{i=1}^{n}=\cup_{k=1}^{q}G_{i_{k}},
 $$
 where $1\leq i_{k}\leq n$ and $1\leq q\leq n.$ Each element in $G_{i_{k}}$ extends analytically to the other
 elements
  in
 $G_{i_{k}}$, but it does not extend to any element in
 $G_{i_{l}}$ if $i_{k}\neq i_{l}$.
Suppose  $\Gamma$ is the curve constructed in Section \ref{Riems},
drawn through these branch points and a fixed point on the unit
circle so that $\D\backslash\Gamma$ is a simply connected region.
For each
$1\leq k\leq q$, define
 a bounded linear operator ${\mathcal E}_{k}: L^{2}_{a}\rightarrow
 L^{2}_{a}$ by
 \begin{equation}\label{Erk}
 ({\mathcal E}_{k}f)(z)=\sum_{\rho\in G_{i_{k}}}\rho^{\prime}(z)f(\rho (z)) \end{equation}
 for $z\in E$ and each $f\in L_{a}^{2}.$   As in the proof of Theorem \ref{Riem}, the operator ${\mathcal E}_{k}$ is
 well-defined since the $\rho (z)$ can be extended to all of $E$ and
 $f(\rho (z))$ is locally bounded on $E$ and hence analytic on $E$.

We observe that the label of each element in $G_{i_{k}}$ depends on
$\Gamma$ and a neighborhood of $z_{0}$ but the set   $G_{i_{k}}$
does not. So we can view   $G_{i_{k}}$ as a collection of local
inverses defined locally and each element in $G_{i_{k}}$ may have
different labels as a global function on $E/\Gamma$. Thus the
summation in (\ref{Erk}) defining  ${\mathcal E}_{k}f(z)$    does
not depend on how we choose the curve $\Gamma$ or how we label those
local inverses on a neighborhood of any point $z_{0}$ in $E$.

First we introduce the notation $\hat{\Gamma}$ which denotes the set
$$\Gamma \cup \cup_{j=1}^{n}\{w\in E\backslash \Gamma :
\rho_{j}(w)\in \Gamma \},$$ which consists of finitely many curves
on $\D$. Since $\rho$ maps $E$ into $E$ and is locally analytic and
injective on $E$, $\hat{\Gamma}$ is a closed set which has area
measure equal to zero.

\begin{lem}\label{ontorho} For each $\rho\in\{\rho_{j}\}_{j=1}^{n},$
$\rho (E\backslash \Gamma)$ contains $\D\backslash \hat{\Gamma}$,
and hence the area measure of $\D\backslash\rho (E\backslash \Gamma
)$ equals zero.
\end{lem}
\proof Let  $z_{0}$ be a point in $\D\backslash \hat{\Gamma}$. Thus
$\rho_{j} (z_{0})$ is not in $\Gamma$ for each $\rho_{j}$. Since
$\rho_{j}$ is locally analytic and injective in $E$, $\rho_{j}$ is
an open mapping. Thus we can find a small open neighborhood
${\mathcal V}_{j}$ of the point $z_{0}$ such that
$$\rho_{j}({\mathcal
V}_{j})\cap \Gamma =\emptyset .$$ Let
$${\mathcal V}=\cap_{j=1}^{n}{\mathcal V}_{j}.$$
Then $\mathcal V$ is an open neighborhood of $z_{0}$ such that
$$\cup_{j=1}^{n}\rho_{j}({\mathcal V})\subset \cup_{j=1}^{n}\rho_{j}({\mathcal V}_{j})\subset E\backslash \Gamma .$$
For each $\rho$, there is a $\hat{\rho}\in \{\rho_{j}\}_{j=1}^{n}$
such that $$\rho(\hat{\rho}(z))=\rho_{1}(z)=z$$ on an open
neighborhood $\hat{\mathcal V}$ of $z_{0}$ which is contained in $
{\mathcal V}$. Letting $w_{0}=\hat{\rho}(z_{0})$, then $w_{0}$ is
contained in $\hat{\rho}({\mathcal V})$ and hence contained in
$E\backslash \Gamma .$ Therefore $z_{0}$ is contained in $\rho
(E\backslash \Gamma )$. This completes the proof.

Since each local inverse $\rho$ is locally injective, there
 are a family of disjoint open sets in $\mathbb D$ on which
 each    $\rho$ in $\{\rho_{j}\}_{j=1}^{n}$ is injective on $U_{\mu}$ and
$$|{\mathbb D}/\cup_\mu{U_{\mu}}|=0.$$
Let $S_{k}$ be the connected component of the Riemann surface
$S_{\phi}$ associated with $G_{i_{k}}$.

\begin{lem}\label{Inje} Let $\{U_{\mu}\}$ be disjoint open sets in $\mathbb D$ on which
 each    $\rho$ in $\{\rho_{j}\}_{j=1}^{n}$ is injective on $U_{\mu}$ and
$$|{\mathbb D}/\cup_\mu{U_{\mu}}|=0.$$  Then
$$|\rho (U_{\mu_{1}})\cap \rho (U_{\mu_{2}})|=0$$
if $\mu_{1}\neq \mu_{2} $.
\end{lem}
\proof As we pointed out in the proof of Theorem \ref{Riem}, each
local inverse $\rho$ has analytic continuation in $r{\mathbb
D}\backslash {\mathcal F}$ for some $r>1$ and $|\rho (z)|=1$ on the
unit circle. So we can construct a Riemann surface ${\mathcal
S}_{r}$ over $r{\mathbb D}\backslash {\mathcal F}$ and $S_{k}$ is
an open region of ${\mathcal S}_{r}$. Assume $G_{i_{k}}$ consists of
$n_{k}$ elements $\tilde{\rho}_{1}, \cdots , \tilde{\rho}_{n_{k}}$. For
 each open set $U$ contained in the unit disk, on which
$\tilde{\rho}_{j}$ is injective,  we define a  function $f$ on
$S_{k}$ by
$$f(\tilde{\rho}_{j}(z), z)=\tilde{\rho}_{j}(z)$$
for each $(\tilde{\rho}_{j}(z), z)\in S_{k}$.   Then $f$ is a
holomorphic function on $S_{k}$. Clearly, $f$ extends to be a
holomorphic function on ${\mathcal S}_{r}.$ Let $\omega$ be the
differential $1$-form on $S_{k}$
$$\omega =\frac{-i}{2}\overline{f}df .$$
On the chart ${\mathcal U}_{\rho}=\{(\rho (z), z):z\in U\}$, it is
easy to check
$$\omega =\frac{-i}{2}\overline{\rho}d\rho $$
and \begin{eqnarray*}
d\omega &=&\frac{-i}{2}d(\bar{f}df)=d{\bar{f}}\wedge df\\
&=&\frac{-i}{2}(\overline{\rho^{\prime}(z)}d\bar{z})\wedge
(\rho^{\prime}(z)dz)\\
&=&\frac{-i}{2}|\rho^{\prime}(z)|^{2}d\bar{z}\wedge dz.
\end{eqnarray*}

  Now $\{\rho (U_{\mu}):\mu
, \rho\in G_{i_{k}}\}$ forms a local chart of the Riemann surface
$S_{k}$ minus a set with zero surface area. Thus
$$|S_{k}|=\sum_{\rho\in G_{i_{k}}}\int_{\cup_{\mu}U_{\mu}}\frac{i}{2}d{\rho} \wedge d\bar{\rho}
=\sum_{\rho\in G_{i_{k}}}\sum_{\mu}|\rho
(U_{\mu})|$$ and
$$|S_{k}|=\int_{S_{k}}d\omega .$$

On the other hand, for each point $\beta \in \mathcal F$ and  small
positive number $\epsilon$, let $C(\beta, {\epsilon})$ be the circle
with center $\beta$ and radius $\epsilon$ and $\Gamma_{j,
\epsilon}(\beta )=\{(\tilde{\rho}_{j}(z), z):z\in C(\beta,
{\epsilon})\}.$ Let $S_{k, \epsilon}=S_{k}\backslash
\cup_{j=1}^{n_{k}}\cup_{\beta\in\mathcal F}\Gamma_{j,
\epsilon}(\beta ).$ Clearly, $S_{k, \epsilon}$ is an open Riemann
surface with boundary
$$\partial S_{k, \epsilon}=[\cup_{j=1}^{n_{k}}\cup_{\beta\in\mathcal
F}\Gamma_{j, \epsilon}(\beta )]\cup \{|f(w,z)|=1: (w,z)\in S_{k}\}$$
in ${\mathcal S}_{r}.$  By Stokes' formula (\cite{CCL}, Theorem
4.2), we have
\begin{eqnarray*}
|S_{k, \epsilon}|&=& \int_{S_{k, \epsilon}}d\omega \\
                 &=& \int_{\partial S_{k, \epsilon}} \omega \\
                 &=&\sum_{j=1}^{n_{k}}\sum_{\beta\in\mathcal F}\int_{\Gamma_{j,
                 \epsilon}(\beta )}\omega +\int_{\{|f(w,z)|=1: (w,z)\in
                 S_{k}\}}\omega\\
                 &=&
                  -\frac{i}{2}\sum_{j=1}^{n_{k}}\sum_{\beta\in {\mathcal F}}
                  \int_{C(\beta
                  ,\epsilon)}\overline{\tilde{\rho}_{j}}d\tilde{\rho}_{j}-\frac{i}{2}\sum_{j=1}^{n_{k}}\int_{|z|=1}
                  \frac{d\tilde{\rho}_{j}(z)}{\tilde{\rho}_{j}(z)}\\
                  &=&-\frac{i}{2}\sum_{j=1}^{n_{k}}\sum_{\beta\in {\mathcal F}}
                  \int_{C(\beta
                  ,\epsilon)}\overline{\tilde{\rho}_{j}}d\tilde{\rho}_{j}-\frac{i}{2}\int_{|z|=1}d[
                  \ln
                  (\prod_{j=1}^{n_{k}}\tilde{\rho}_{j}(z))].
 \end{eqnarray*}
Noting that the product $\prod_{j=1}^{n_{k}}\tilde{\rho}_{j}(z)$ is
a Blaschke factor of $\frac{\phi (0)-\phi}{1-\overline{\phi
(0)}\phi}$ with order $n_{k}$, by the formula of the winding number
of a closed curve, we have
$$-\frac{i}{2}\int_{|z|=1}d[
                  \ln
                  (\prod_{j=1}^{n_{k}}\tilde{\rho}_{j}(z))]=n_{k}\pi
                  =n_k |{\mathbb D}|
                  .$$
                  By the Puiseux theorem (\cite{Bli}, Lemma 13.1 and Theorem 13.1 or \cite{Fos}, Theorem 8.14), noting that each local
                  inverse $w=\rho (z)$ is a solution
                  $$P(w)Q(z)-P(z)Q(w)=0,$$
                  and the leading coefficient of $w^{n}$ in the above polynomial is given by
                  $$Q(z)-P(z)\overline{\phi (0)},$$
                  which never vanishes on the unit disk,  for each $\beta\in
                  \mathcal F$. Therefore  we see that there is a neighborhood $D(\beta ,
                  \epsilon )\backslash \{\beta \}$, at which $\rho (z)$ has a
                  power series expansion of $(z-\beta )^{1/n_{k}}.$
                  Thus for some positive constant $M>0$,
$$|\int_{C(\beta
                  ,\epsilon)}\overline{\tilde{\rho}_{j}}d\tilde{\rho}_{j}|\leq
                  M\epsilon^{1/n_{k}}\rightarrow 0,$$
                  and
                  $$|\int_{D(\beta ,
                  \epsilon)\backslash \{\beta\}}d\tilde{\rho}_{j}(z)\wedge
                  d\overline{\tilde{\rho}}_{j}(z)|\leq
                  M\epsilon^{2/n_{k}}\rightarrow 0,$$
  as $\epsilon\rightarrow 0$.
                  This gives
                  $$\lim_{\epsilon\rightarrow 0}|S_{k,
                  \epsilon}|=n_{k}|{\mathbb D}|,$$
                  and
                  $$\lim_{\epsilon\rightarrow}|S_{k,
                  \epsilon}|=|S_{k}|.$$
    Hence
$$|S_{k}|=n_{k}|{\mathbb D}|.$$
 By Lemma \ref{ontorho},
we have
$$|\rho ({\mathbb D}/\Gamma )|=|{\mathbb D}|.$$
This gives
$$\sum_{\mu}|\rho
(U_{\mu})|\geq |\rho ({\mathbb D}/\Gamma )|=|{\mathbb D}|.$$ Thus
$$|S_{k}|=\sum_{\rho\in G_{i_{k}}}\sum_{\mu}|\rho
(U_{\mu})|\geq n_{k}|{\mathbb D}|.$$  But this implies
$$\sum_{\mu}|\rho
(U_{\mu})|=|{\mathbb D}|.$$ Therefore
$$|\rho (U_{\mu_{1}})\cap \rho (U_{\mu_{2}})|=0$$
for $\mu_{1}\neq \mu_{2}.$ This completes the proof.

 The boundedness
of ${\mathcal E}_{k}$ follows from
 \begin{eqnarray}\label{change}
 \int_{E\backslash \Gamma}|f(\rho
 (z))|^{2}|\rho^{\prime}(z)|^{2}dA(z) &=&\sum_{\mu}\int_{\rho
 (U_{\mu})}|f(w)|^{2}dA(w)\notag\\
&=&\int_{{\mathbb D}}|f(w)|^{2}dA(w),
 \end{eqnarray}
where the $\{U_{\mu}\}$ are disjoint open sets in $\mathbb D$ on which
each    $\rho$ is injective on $U_{\mu}$ and
$$|{\mathbb D}/\cup_\mu{U_{\mu}}|=0.$$  The last equality in (\ref{change}) follows from the above lemma.

 To get ${\mathcal E}_{k}^{*}$,
we need the following change of variable formula.

\begin{lem}
For each $\rho\in\{\rho_{j}\}_{j=1}^{n}$ and $f\in
 L^{2}_{a},$
 \begin{eqnarray*}
 \int_{E\backslash \Gamma}|f(\rho
 (z))|^{2}|\rho^{\prime}(z)|^{2}dA(z)
&=& \int_{\D}|f(w)|^{2}dA(w).
 \end{eqnarray*}
\end{lem}
\proof We can choose $\{U_{\mu}\}$ to be  disjoint open sets in
$\mathbb D$ such that for each $\rho$,  $\rho$ is injective on
$U_{\mu}$ and
$$|{\mathbb D}/\cup_\mu {U_{\mu}}|=0.$$ For each $\rho\in\{\rho_{j}\}_{j=1}^{n}$ and $f\in
 L^{2}_{a}$, we have
\begin{eqnarray*}
 \int_{E\backslash \Gamma}|f(\rho
 (z))|^{2}|\rho^{\prime}(z)|^{2}dA(z) &=&\int_{E\backslash \Gamma}|f(\rho
 (z))|^{2}|\rho^{\prime}(z)|^{2}dA(z)\\
 &=&\sum_{\mu}\int_{\rho (U_{\mu})}|f(w)|^{2}dA(w)\\
 &=&\int_{\rho (E\backslash \Gamma)}|f(w)|^{2}dA(w)\\
&=& \int_{\D}|f(w)|^{2}dA(w).
 \end{eqnarray*}
The first equality follows from the fact that the area measure of
$\D\backslash E$ is zero. The second equality follows from Lemma
\ref{Inje}.  The last equality comes from   Lemma
\ref{ontorho}, which states that the area measure of $\D\backslash
\rho (E\backslash \Gamma)$ is zero. This completes the proof.

Note that $\rho_{1}$ equals $z$. Since $\rho^{-1}$ is also in $
\{\rho_{i}\}_{i=1}^{n}$ for each $\rho\in \{\rho_{i}\}_{i=1}^{n},$
let $G_{i_{k}}^{-}$ denote the subset of $\{\rho_{i}\}_{i=1}^{n}$:
$$G_{i_{k}}^{-}=\{\rho :\rho^{-1}\in G_{i_{k}}\}.$$
We need the following lemma to find ${\mathcal E}_{k}^{*}$.

\begin{lem}\label{Ginverse} For each $i_{k}$, there is an integer $k^{-}$ with $1\leq k^{-}\leq q$ such that
$$G_{i_{k}}^{-}=G_{i_{k^{-}}}.$$
\end{lem}

\proof Assume that $z_{0}$ is a point in $E\backslash \Gamma .$ For
two elements $\rho$ and $\hat{\rho}$ in $G_{i_{k}},$ suppose that
$\hat{\rho}$ is an analytic continuation of $\rho$ along some loop
$\gamma$  in $E$ containing $z_{0}$. Let $\hat{\gamma}$ be the image
curve of $\gamma$ under this analytic continuation along $\gamma$.
We will show that $\hat{\rho}^{-1}$ is an analytic continuation of
$\rho^{-1}$ along a loop at $z_{0}$ in $E$.

  Note that $\hat{\gamma}$ is a curve connecting $\rho (z_{0})$ and $\hat{\rho}(z_{0})$ in $E$.
Thus $\rho^{-1}$ has an analytic continuation $\tilde{\rho}$ from $\rho (z_{0})$ to $\hat{\rho}(z_{0})$ along
 the curve $\hat{\gamma}$ in $E$.
This gives that $\rho_{1}=\rho^{-1}\circ\rho $ has an analytic
continuation $\tilde{\rho}\circ\hat{\rho}$ along the loop $\gamma$
from $z_{0}$ to $z_{0}$ in $E$, but $\rho_{1}=z$ has only one
analytic continuation $\rho_{1}$ along any loop in $E$. Hence we
have that $\tilde{\rho}\circ\hat{\rho}=\rho_{1}$ to get
$\tilde{\rho}=\hat{\rho}^{-1}$. This means that $\rho^{-1}$ can be
analytically extended to $\hat{\rho}^{-1}$ along a curve connecting
$\rho (z_{0})$
 and
$\hat{\rho}(z_{0})$ in $E$. Also we can find two curves $\gamma_{1}, \gamma_{2}$ in $E\backslash\Gamma$ such that
$\rho^{-1}$ on an open neighborhood of $z_{0}$ is an  analytic continuation of $\rho^{-1}$
on an open neighborhood of $\rho (z_{0})$ along $\gamma_{1}$ and $\hat{\rho}^{-1}$ on an open neighborhood of
$z_{0}$ is an  analytic continuation of $\hat{\rho}^{-1}$
on an open neighborhood of $\hat{\rho} (z_{0})$ along $\gamma_{2}$ . Thus $\hat{\rho}^{-1}$ is an analytic continuation
of $\rho^{-1}$ along the loop $\gamma_{1}\cup\hat{\gamma}\cup (-\gamma_{2})$ in $E$.

As $\{\rho_{i}\}_{i=1}^{n}$ has a group-like property under the
composition defined above, there is an integer $k^{-}$ with $1\leq
k^{-}\leq q$ such that $G_{i_{k}}^{-}$ equals  $G_{i_{k^{-}}}$. This
completes the proof.

\begin{lem}\label{ek} For  each integer $k$ with $1\leq k \leq q$, there is an integer
$k^{-}$ with $1\leq k^{-}\leq q$ such that $${\mathcal
E}_{k}^{*}={\mathcal E}_{k^{-}}.$$
\end{lem}
\proof Polarizing  the   change  of variable formula (\ref{change})
gives
 $$\int_{E\backslash \Gamma}f(\rho (z))\overline{g(\rho
 (z))}|\rho^{\prime}(z)|^{2}dA(z)=
 \int_{\D}f(w)\overline{g(w)}dA(w).$$
for two polynomials $f, g$ of $z$.  Choose a collection
$\{U_{\mu}\}$ of disjoint open subsets such that $$ |E/[\Gamma
\cup_{\mu}U_{\mu}]|=0,$$  for each $\mu$,  the sets $\{\rho
(U_{\mu})\}_{\rho\in G_{i_{k}}}$  are disjoint and for each $\rho$
and $\mu$, $\rho$ is injective on $U_{\mu}$.
 Note that for
each point $\rho (z)\notin \rho (\Gamma )$, by Lemma \ref{Ginverse},
there is a $\hat{\rho}_{\mu}\in G_{i_{k^{-}}}$ such that
$$\hat{\rho}_{\mu}(\rho (z))=z$$
on a neighborhood $U_{\mu}$ of $z$. Let $w=\rho (z)$. Since $\rho$
is analytic  and locally injective on $E$, the above equality gives
that there is a $\tilde{\rho}_{\mu}\in G_{i_{k}}$ such that
$$\tilde{\rho}_{\mu} (\hat{\rho}_{\mu} (w))=w$$
for $w\in \rho (U_{\mu})$. Hence
$$\rho (z)= \tilde{\rho}_{\mu}(z)$$
for $z\in U_{\mu}$ and
$$\rho ( \hat{\rho}_{\mu} (w))=w$$
on a neighborhood $\rho (U_{\mu})=\tilde{\rho}_{\mu}(U_{\mu})$ of
$\rho (z_{0})$. Let $\chi_{\rho (U_{\mu})}$ denote the
characteristic function of the set $\rho (U_{\mu}).$ Thus
\begin{eqnarray*}
\ip{{\mathcal E}_{k}^{*}g,f} &=& \ip{g, {\mathcal E}_{k} f}\\
                             &=& \int_{\D}\sum_{\rho\in
                             G_{i_{k}}}g(z)\overline{f(\rho
                             (z))\rho^{\prime}(z)}dA(z)\\
                             &=& \int_{E\backslash \hat{\Gamma}}
                            \sum_{\rho\in
                             G_{i_{k}}}g(z)\overline{f(\rho
                             (z))\rho^{\prime}(z)}dA(z)\\
                             &=& \sum_{\rho\in
                             G_{i_{k}}}\sum_{\mu}\int_{U_{\mu}}
                              g(\hat{\rho}_{\mu}\circ \rho (z))\overline{f(\rho
                             (z))\rho^{\prime}(z)}dA(z)\\
                             &=&
                            \sum_{\rho\in G_{i_{k}}}\sum_{\mu}\int_{U_{\mu}}
                           \frac{g(\hat{\rho}_{\mu}\circ \rho (z))}{\rho^{\prime}(z)}
                             \overline{f(\rho
                             (z)) }|\rho^{\prime}(z)|^{2}  dA(z)\\
                              &=& \sum_{\rho\in G_{i_{k}}}\sum_{\mu}\int_{\rho (U_{\mu})}
                              \frac{g(\hat{\rho}_{\mu}(w))
                            }{\rho^{\prime}(\hat{\rho}_{\mu}(w))}
                             \overline{f(w) }  dA(w)\\
                             &=&\int_{{\mathbb D}}\sum_{\rho\in G_{i_{k}}}\sum_{\mu}\chi_{\rho
                             (U_{\mu})}(w)
 {g(\hat{\rho}_{\mu}(w))
                            }{\hat{\rho}_{\alpha}^{\prime}(w)}
                             \overline{f(w) }  dA(w).
\end{eqnarray*}
The third equality follows from the fact that $\hat{\Gamma}$ has
  area measure equal to zero and the fact that  $\sum_{\rho\in
                             G_{i_{k}}}f(\rho
                             (z))\rho^{\prime}(z)$ is in the Bergman
                             space $L^{2}_{a}$. The last equality follows from
 $$\hat{\rho}_{\mu}^{\prime}(w)=\frac{1}{\rho^{\prime}(\hat{\rho}_{\mu}(w))}.$$
 Let $S_{k}$ be the connected component of the Riemann surface for $\phi^{-1}\circ\phi$ over $\mathbb D$ associated with
 $G_{i_{k}}$ and $n_{k}$ the cardinality of $G_{i_{k}}$.  Since $S_{k}$ is an $n_{k}$-sheeted ramified covering of
    $\mathbb D$, it can be pictured as $n_{k}$ unit disks attached with appropriate branch points and lying over $\mathbb D$.
     By Lemma \ref{Inje},
   for almost each $w$ in $\mathbb D$, there are only $n_{k}~~$ elements $\mu_{i}$ such that
$$w=\rho_{\mu_{i}}(z_{\mu_{i}})$$
for each $z_{\mu_{i}}\in U_{\mu_{i}}$, $i=1, \cdots , n_{k}.$ Thus,
for almost  all $w\in {\mathbb D}$,
\begin{eqnarray*}
\sum_{\rho\in G_{i_{k}}}\sum_{\mu}\chi_{\rho (U_{\mu})}(w)
 {g(\hat{\rho}_{\mu}(w))
                            }{\hat{\rho}_{\mu}^{\prime}(w)}&=&\sum_{i=1}^{n_{k}}
 {g(\hat{\rho}_{\mu_{i}}(w))
                            }{\hat{\rho}_{\mu_{i}}^{\prime}(w)}.
                            \end{eqnarray*}

We claim that for each $w$, $\hat{\rho}_{\mu_{i}}(w)\neq
\hat{\rho}_{\mu_{j}}(w)$ if $i\neq j$. If this is not true, for some
$i\neq j$,
$$\hat{\rho}_{\mu_{i}}(w)=\hat{\rho}_{\mu_{j}}(w),$$
as
$$w=\rho_{\mu_{i}}(z_{\mu_{i}})=\rho_{\mu_{j}}(z_{\mu_{j}}).$$
Thus we have \begin{eqnarray*}
z_{\mu_{i}}&=&\hat{\rho}_{\mu_{i}}(\rho_{\mu_{i}}(z_{\mu_{i}}))\\
              &=&\hat{\rho}_{\mu_{i}}(w)\\
              &=&\hat{\rho}_{\mu_{j}}(w)\\
              &=&\hat{\rho}_{\mu_{j}}(\rho_{\mu_{j}}(z_{\mu_{j}}))\\
              &=& z_{\mu_{j}}.
              \end{eqnarray*}
This implies that
$\rho_{\mu_{i}}(z_{\mu_{i}})=\rho_{\mu_{j}}(z_{\mu_{i}})$, which
contradicts the fact that the intersection of  $\rho_{\mu_{i}}(U_{\mu_{i}})$ and $ \rho_{\mu_{j}}(U_{\mu_{i}})$ is empty.  Thus the
$\rho_{\mu_{j}}(U_{\mu_{i}})$ are disjoint.

Since $G_{i_{k^{-}}}$ has $n_{k}$ elements and contains
$\{\hat{\rho}_{\mu_{i}}\}_{i=1}^{n_{k}},$ we have
$$\sum_{i=1}^{n_{k}}
 {g(\hat{\rho}_{\mu_{i}}(w))
                            }{\hat{\rho}_{\mu_{i}}^{\prime}(w)}=\sum_{\rho\in G_{i_{k^{-}}}}
 {g(\rho(w))
                            }{\rho^{\prime}(w)} $$
even if those $\mu_{i}$ may depend on the point $w$.

    If $g$ is a polynomial of $z$, by
the proof of Theorem \ref{Riem}, the right hand side of the above
equality is analytic in $E$. Letting $V$ be an open neighborhood of
a branch point of $\phi$, we have
$$\int_{V\backslash \Gamma}|\rho^{\prime}(z)|^{2}dA(z)=|\rho (V\backslash \Gamma )|<\infty .$$
This implies that $\sum_{\rho\in G_{i_{k^{-}}} }\rho^{\prime}
 (z)g(\rho (z))$   extends analytically on $\D$ and is
 in the Bergman space $L^{2}_{a}.$
 Thus
 \begin{eqnarray*}
 \ip{{\mathcal E}_{k}^{*}g,f} &=&\int_{{\mathbb D}}\sum_{\rho\in G_{i_{k}}}\sum_{\mu}\chi_{\rho (U_{\mu})}
 {g(\hat{\rho}_{\mu}(w))
                            }{\hat{\rho}_{\mu}^{\prime}(w)}
                             \overline{f(w) }  dA(w)\\
                             &=&\int_{\D}
                             [\sum_{\rho\in G_{i_{k^{-}}}}
g({\rho}(w)){\rho}^{\prime}(w)]
                             \overline{f(w) }  dA(w)\\
                             &=& \ip{{\mathcal E}_{k^{-}}g,f}.\\
                             \end{eqnarray*}
  Since  $\sum_{\rho\in G_{i_{k^{-}}}}
g({\rho}(w)){\rho}^{\prime}(w)$ is in the Bergman space and  the
polynomials are dense in the Bergman space $L^{2}_{a}$, we have that
for any polynomial $g$,
\begin{eqnarray*}
{\mathcal E}_{k}^{*}g(z) &=& \sum_{\rho\in
G_{i_{k^{-}}}}\rho^{\prime}
 (z)g(\rho (z)),
 \end{eqnarray*}
and hence
\begin{eqnarray*}
{\mathcal E}_{k}^{*}g(z)
 &=&{\mathcal E}_{k^{-}}g(z)
 \end{eqnarray*}
 for $z\in \D$. By the fact that the polynomials are dense in the Bergman space $L^{2}_{a}$,  we have that ${\mathcal E}_{k}^{*}=
 {\mathcal E}_{k^{-}}.$ This completes the proof.

 \begin{thm}\label{mthlast} Let $\phi$ be a finite Blaschke product. The von Neumann algebra ${\mathcal A}_{\phi}$ is
 generated by the linearly independent operators ${\mathcal E}_{1}, \cdots , {\mathcal E}_{q}$ and hence has dimension $q$.
\end{thm}

\proof Let $q$ be the number of connected components of the Riemann
surface $\phi^{-1}\circ\phi$ over the unit disk. Recall that
${\mathcal A}_{\phi}$ is the von Neumann algebra
$\{M_{\phi}\}^{\prime}\cap \{M_{\phi}^{*}\}^{\prime}.$

To finish the proof we need show that ${\mathcal A}_{\phi}$ is a
finite dimensional space with dimension equal to $q$.

  By Lemma \ref{decomp},  for each unitary operator $W$ in ${\mathcal A}_{\phi}$,
there are at most $q$ distinct complex numbers $\hat{r}_{1}, \cdots
, \hat{r}_{q}$ such that for each $\alpha$ in $U$,
  \begin{eqnarray*}
    W^{*}k_{\alpha}&=& \sum_{k=1}^{q}\hat{r}_{k}\sum_{\rho\in G_{i_{k}}}\overline{\rho^{\prime}(\alpha )}k_{\rho (\alpha )},
  \end{eqnarray*}
  where $\hat{r}_{k}=r_{\rho}$ for some $\rho$ in $G_{i_{k}}.$  Since $\{k_{\alpha}\}_{\alpha\in U}$ is dense in the Bergman space,
  we have
  $$W^{*}=\sum_{k=1}^{q}\hat{r}_{k}{\mathcal E}_{k}.$$ Thus
  ${\mathcal A}_{\phi}$ contains at most $q$ linearly independent unitary
  operators. By the Russo-Dye Theorem, \cite{RuD}, \cite{Dav},
\cite{Zh2},  every element in ${\mathcal A}_{\phi}$ can be written
as a finite linear combination of unitary operators in ${\mathcal
U}_{{\mathcal A}_{\phi}}$. Thus ${\mathcal A}_{\phi}$ is a finite
dimensional space with dimension at most $q$.

 Next we  show that the dimension of ${\mathcal A}_{\phi}$
 is at least  $q.$
By (\ref{Erk}), there are $q$ bounded linear operators ${\mathcal
E}_{1}, \cdots , {\mathcal E}_{q}$ on the Bergman space.
     Since $\phi (\rho (z))=\phi
 (z)$, we have that
 $$M_{\phi}{\mathcal E}_{k}={\mathcal E}_{k}M_{\phi}.$$
    Thus the ${\mathcal E}_{1}, \cdots , {\mathcal E}_{q}$ are contained in
    $\{M_{\phi}\}^{\prime}.$ Now
    Lemma \ref{ek} tells us that the ${\mathcal E}_{1}^{*}, \cdots , {\mathcal E}_{q}^{*}$ are also contained in $\{M_{\phi}\}^{\prime}$.
    This gives that  the ${\mathcal E}_{1}, \cdots , {\mathcal E}_{q}$ are contained in
    ${\mathcal A}_{\phi}.$  We claim
 that the ${\mathcal E}_{1}, {\mathcal E}_{2}, \cdots , {\mathcal E}_{q}$
are linearly
 independent. If this is false, there are constants $c_{1}, \cdots ,
 c_{q}$, not all of which   are zero,  such that
 $$c_{1}{\mathcal E}_{1}+\cdots + c_{q}{\mathcal E}_{q}=0.$$
 Thus for each $\alpha$ in $E$, we have
 $$[c_{1}{\mathcal E}_{1}+\cdots + c_{q}{\mathcal
 E}_{q}]^{*}k_{\alpha}=0.$$
On the other hand, ({\ref{Erk}) gives
$$[c_{1}{\mathcal E}_{1}+\cdots + c_{q}{\mathcal
 E}_{q}]^{*}k_{\alpha}=\sum_{k=1}^{q}\overline{c}_{k}\sum_{\rho\in
 G_{i_{k}}}\overline{\rho^{\prime} (\alpha)}k_{\rho (\alpha)}.$$
 As in the proof of Theorem \ref{localrepu}, for each $i$, define
  $$P_{i}(\alpha , z)=\prod_{j\neq i}^{n}(z-\rho_{j}(\alpha )).$$
  An easy calculation gives
$$\ip{P_{i_{k}}(\alpha , .),[c_{1}{\mathcal E}_{1}+\cdots + c_{q}{\mathcal
 E}_{q}]^{*}k_{\alpha}}=c_{k}\rho_{i_k}^{\prime}(\alpha )P_{i_{k}}(\alpha, \rho_{i_{k}}(\alpha
 )).$$
 Since  $P_{i_{k}}(\alpha, \rho_{i_{k}}(\alpha
 ))\neq 0$ and $\rho^{\prime}_{i_{k}}(z)$ vanishes only on a countable subset of $\D$,
 we have that  $c_{k}$ must be zero for each $k$. This is a
 contradiction. We conclude that the  ${\mathcal E}_{1}, {\mathcal E}_{2}, \cdots , {\mathcal E}_{q}$
are linearly
 independent to obtain that the dimension of ${\mathcal A}_{\phi}$
 is at least  $q.$
 This completes the proof.

\section{Abelian  von Neumann algebra ${\mathcal A}_{\phi}$}

 In this section, we will show that the von Neumann algebra ${\mathcal A}_{\phi}$ is abelian if the order of
 the Blaschke product $\phi$ is smaller than or equal to $8$. First we recall some concepts and notation from previous sections.  For each $z\in E$, the function $\phi$ is one-to-one in some open neighborhood $D_{z_{i}}$ of each point $z_{i}$ in  $\phi^{-1}\circ\phi (z)=\{z_{1}, \cdots , z_{n}\}$.  Let  $\phi^{-1}\circ\phi=\{\rho_{k}(z)\}_{k=1}^{n}$ be $n$ solutions of the equation $\phi (\rho (z))=\phi (z).$ Then $\rho_{j}(z)$ is
 locally analytic and arbitrarily continuable in $E$. Assume that $\rho_{1}(z)=z$. Every open subset $V$ of $E$ is
 invertible for $\phi$.
 Then $\{\rho_{j}\}_{j=1}^{n}$ is the family of  admissible local inverses in some invertible open disc $V\subset \mathbb D$.
 For a given point $z_{0}\in V$, label the local inverses as $\{\rho_{j}(z)\}_{i=1}^{n}$ on $V$.
 If there is a loop $\gamma$ in $ E$ at $z_{0}$ such that $\rho_{j}$ and  $\rho_{j^{\prime}}$ in $\{\rho_{i}(z)\}_{i=1}^{n}$
 are mutually analytically continuable along $\gamma$, we   write
$$\rho_{j}\thicksim \rho_{j^{\prime}},$$
and it is easy to check that $\thicksim$ is an equivalence relation.
Let $\phi$ be a finite Blaschke product. Let
$ {G}$ be those local inverses of $\phi$ which extend analytically
to only themselves in ${\mathbb D}\backslash {\mathcal F}.$

\begin{lem} $
{G}$ is an elementary subgroup of $Aut({\mathbb D})$ consisting of
elliptic M\"obius transforms and identity $\rho_{1}$. Moreover,  there are
a point $\alpha$ in $\mathbb D$, a unimodular constant $\lambda$ and an integer $n_{G}$ such that
$$G=\{\frac{|\alpha |^{2} -\lambda}{1-|\alpha
|^{2}\lambda}\phi_{\frac{ \alpha (1-\bar{\lambda})}{1-|\alpha
|^{2}\bar{\lambda}}}(z):\lambda^{n_{G}}=1\}.$$
\end{lem}

\proof Since  $G$ consists of these local inverses  $\rho$  which
are equivalent only to themselves, each $\rho$ in $G$ has
analytic continuation on $E$.  Therefore each $\rho$ has an analytic continuation on the
unit disk. Also $|\rho (z)|=1$ on the unit circle and thus $\rho$ is
an inner function. On the other hand, $\rho$ is locally injective. We
conclude that each $\rho$ in $G$ is in $Aut ({\mathbb D})$. For two
elements $\rho$ and $\tau$ in $G$, $\rho\circ\tau$ is still a local
inverse in $G$, which shows that $G$ is a finite subgroup of $Aut
({\mathbb D})$. As $G$ is a finite group, letting $n_{G}$ be the
number of elements in $G$, for each $\rho$ in $G$, we have
$$\rho^{\circ n_{G}}=\rho_{1}(z).$$
Thus $\rho$ is elliptic and so $G$ is elementary.

 According to (\cite{Brd} on page 12), $G$ has an invariant point $\alpha$ in
$\mathbb D$ and so it is a group of hyperbolic rotations about
$\alpha$, that is
$$G=\{\phi_{\alpha}(\lambda \phi_{\alpha} (z)):
\lambda^{n_{G}}=1\}=\{\frac{|\alpha |^{2} -\lambda}{1-|\alpha
|^{2}\lambda}\phi_{\frac{ \alpha (1-\bar{\lambda})}{1-|\alpha
|^{2}\bar{\lambda}}}(z):\lambda^{n_{G}}=1\}.$$
This completes the
proof.

{\bf Remark.} Clearly, the group $G$ in the above lemma is abelian.

\begin{lem}\label{dist}
Suppose that $\lambda$ is the $k$-th root of  unity and $\alpha$
is a nonzero point in the unit disk. If $\beta_{j}=\frac{ \alpha
(1-\bar{\lambda}^{j})}{1-|\alpha |^{2}\bar{\lambda}^{j}},$ then
$\{\beta_{1}, \beta_{2}, \cdots , \beta_{k}\}$ are distinct points
in the unit disk.
\end{lem}

\proof  Suppose that $\beta_{j}=\beta_{l}$. Then
$$\frac{ \alpha
(1-\bar{\lambda}^{j})}{1-|\alpha |^{2}\bar{\lambda}^{j}}=\frac{
\alpha (1-\bar{\lambda}^{l})}{1-|\alpha |^{2}\bar{\lambda}^{l}}.$$
Since $\alpha$ does not equal zero, we have
$$\frac{
(1-\bar{\lambda}^{j})}{1-|\alpha |^{2}\bar{\lambda}^{j}}=\frac{
  (1-\bar{\lambda}^{l})}{1-|\alpha |^{2}\bar{\lambda}^{l}},$$
which yields
  $$\bar{\lambda}^{l}-\bar{\lambda}^{j}=|\alpha |^{2}
  (\bar{\lambda}^{l}-\bar{\lambda}^{j}).$$
  Since $\alpha$ is in the open unit disk, the above equality gives
$$\bar{\lambda}^{l}-\bar{\lambda}^{j}=0.$$
If $1\leq j, l\leq k$, then $j$ must equal $l$. This implies that
$\{\beta_{1}, \beta_{2}, \cdots , \beta_{k}\}$ are distinct points.
This completes the proof.

Let $q$ denote the  number of connected components of the Riemann surface
$S_{\phi}.$
Let $\phi$ be a Blaschke product. We say that $\phi$ is reducible if
there are two other  Blaschke products, $\phi_{1}$ and $\phi_{2}$, with
orders larger than $1$ such that $\phi=\phi_{1}\circ\phi_{2}.$

\begin{lem}\label{factor} Let $\phi$ be a Blaschke product with order $n$. If ~~$G$
has $n_{G}>1$ elements, then $\phi$ is reducible and $n_{G}|n$.
\end{lem}

\proof  Let $\lambda$ denote $\phi (0)$ and
$$\psi=\phi_{\lambda}\circ\phi .$$ Since for each $\rho\in G$,
$\phi\circ\rho =\phi$, we have that $\psi\circ\rho =\psi .$ Noting
that $\psi (0)=0$, we write
$$\psi =z\psi_{1}.$$
Since $\psi \circ\rho =\psi$,  for each $\rho\in G$, $\rho$ is a
factor of $\psi$. Letting $\phi_{G}=\prod_{\rho\in G}\rho (z),$ by
Lemma \ref{dist},  we have
$$\psi =[\prod_{\rho\in G}\rho (z)]\tilde{\psi}(z)=\phi_{G}\tilde{\psi}$$
for a Blashcke product $\tilde{\psi}$ with order $n-n_{G},$  to yield
$$\tilde{\psi}\circ\rho =\tilde{\psi}$$
for each $\rho$ in $G$.  Thus
$$\phi=\phi_{\lambda}\circ\psi =\phi_{\lambda}\circ
(\phi_{G}\tilde{\psi}).$$ Repeating the above argument applied to
$\tilde{\psi}$, noting that the order of $\tilde{\psi}$ is
$n-n_{G}$, and  using induction, we obtain
$$\tilde{\psi}=\psi_{2}\circ \phi_{G}$$
for some Blaschke product $\psi_{2}.$ This gives
$$\phi =\phi_{\lambda}\circ (\phi_{G}\psi_{2}\circ
\phi_{G})=\psi_{3}\circ \phi_{G},$$ where
$\psi_{3}=\phi_{\lambda}\circ (z\psi_{2}).$ So $\phi$ is reducible
and the order $n$ of $\phi$ equals the product of  $n_{G}$ and the
order of $\psi_{3}.$ This completes the proof.

\begin{cor} Let $\phi$ be a Blaschke product of order $n$. If $n$ is
greater than or equal to $5$, the number $q$ of connected components
of the Riemann surface $\phi^{-1}\circ\phi$ over the unit disk does
not equal $n-1$
\end{cor}

\proof If $q$ equals $n-1$, then $G$ has $n-2$ elements. Lemma
\ref{dist} gives that $(n-2)|n$, which is impossible.

\begin{thm} Let $\phi$ be a finite Blaschke product with order less
than or equal to $8$. Then ${\mathcal A}_{\phi}$ is commutative and
hence, in these cases, the number of minimal reducing subspaces of
$M_{\phi}$ equals the number of connected components of the Riemann
surface $\phi^{-1}\circ\phi$ over the unit disk.
\end{thm}

\proof As shown in \cite{SZZ1}, the center of the algebra  ${\mathcal
A}_{\phi}$ contains a non trivial projection $P$. So if the
dimension $q$ of the algebra  ${\mathcal A}_{\phi}$ is less than
$5$, by the classification  theorem of finite dimensional von
Neumann algebras \cite{Dav}, ${\mathcal A}_{\phi}$ is commutative. Let $n$ be
the order of the Blaschke product.  So we may assume $n\geq q\geq 5$

If $n=5$ and $q=5$, then the order of $G$ equals $5$ and $G$ is a
cyclic group. Thus $\{{\mathcal E}_{j}\}_{j=1}^{n}$ is commutative
and hence ${\mathcal A}_{\phi}$ is commutative.

If $n=6$ and $q=6$, similarly we have that ${\mathcal A}_{\phi}$ is
commutative.  If $q=5$, $G$ has four elements and hence $n_{G}$
equals $4$.  Lemma \ref{factor} implies that $4$ would be a factor of $6$, which is  impossible.

If $n=7$ and $q=7$, as above, $G$ is a cyclic group and so
${\mathcal A}_{\phi}$ is commutative.  If $q$ equals $6$ or $5$,
then $G$ contains either $5$ or $4$ or $3$ elements. In any of these cases,
Lemma \ref{factor} implies that   $5$, or $4$ or $3$ is a
factor of $7$. This is impossible.

If $n=8$, we consider  four cases $5\leq q\leq  8$.

{\bf Case 1.} $q=8$, then $G$ is a cyclic group and so ${\mathcal
A}_{\phi}$ is commutative.

{\bf Case 2.} $q=7$, then two of local inverses extend analytically each to the other in $E$ and hence $G$ contains $6$ elements. Lemma \ref{factor}
implies that $6$ is a factor of $8$, which is impossible.

{\bf Case 3.}  $q=6$. In this case we will show that the center of
${\mathcal A}_{\phi}$ has dimension at least $4$. So by the
classification  theorem of finite dimensional von Neumann algebras,
${\mathcal A}_{\phi}$ is commutative.

 In this case,  $G$
contains $4$ elements and is a cyclic group generated by an elliptic
M\"obius transform $\rho$. Then the set of local inverses is
divided into $G_{1}=\{\rho_{1}\}$, $G_{2}=\{\rho\}$,
$G_{3}=\{\rho^{2}\}$, $G_{4}=\{\rho^{3}\}$, $G_{5}=\{\rho_{5},
\rho_{6}\}$, and $G_{6}=\{\rho_{7}, \rho_{8}\}$. Let ${\mathcal
E}_{i}$ be the operator associated with $G_{i}$ for $i=1, \cdots ,
6$. Then ${\mathcal E}_{2}$ is a unitary operator and
$${\mathcal E}_{i}={\mathcal E}_{2}^{i-1}$$
for $i=1, \cdots , 4$ and  $${\mathcal E}_{2}^{*}={\mathcal
E}_{2}^{3}.$$ Moreover, ${\mathcal E}_{2}$ commutes with ${\mathcal E}_{i}$
for $i=1,2,3,4$.

For two local inverses, $\tau_{1}$, $\tau_{2}$, if $\tau_{1}$ is
equivalent to $\tau_{2}$, then $\tau_{1}\circ\rho$ is equivalent to
$\tau_{2}\circ\rho$. We observe that both ${\mathcal E}_{i}{\mathcal
E}_{2}$ and ${\mathcal E}_{i}^{*}$ are  in $\{{\mathcal E}_{5},
{\mathcal E}_{6}\}$. Thus there are permutations $\tau$ and $\sigma$
of $\{5, 6\}$ such that
$${\mathcal E}_{i}{\mathcal E}_{2}={\mathcal E}_{\sigma (i)},$$
and $$ {\mathcal E}_{i}^{*}={\mathcal E}_{\tau (i)}.$$ Noting that
$\sigma^{2}(i)=i$, we have
$${\mathcal E}_{i}{\mathcal E}_{2}^{3}={\mathcal E}_{\sigma (i)}.$$
Taking adjoint of both sides of the above equality gives
$$({\mathcal E}_{2}^*)^{3}{\mathcal E}_{\tau (i)}={\mathcal E}_{\tau (\sigma
(i))}.$$ Since $$({\mathcal E}_{2}^*)^{3}={\mathcal E}_{2}$$ and
$\tau$ commutes with $\sigma$, we have
$${\mathcal E}_{2}{\mathcal E}_{ i}={\mathcal E}_{\sigma
(i)}.$$  Thus for $i=5,6,$ $${\mathcal E}_{2}{\mathcal E}_{ i}= {\mathcal E}_{
i}{\mathcal E}_{2}.$$
This gives that ${\mathcal E}_{2}$ is in the
center of ${\mathcal A}_{\phi},$ and so the dimension of the center
is at least $4$.

{\bf Case 4.} $q=5$. Lemma \ref{factor} gives that $G$ has $2$ or
$4$  elements.

If $G$ has $4$ elements, then the local inverses of $\phi$ is
divided into $5$ equivalent classes $G_{1}, \cdots , G_{5}$. We may
assume that $G_{i}=\{\rho^{i-1}(z)\}$ for $i\leq 4,$ where the zero power $\rho^0 (z)$ of $\rho (z)$
is $\rho_{1}(z)$ and $G_{5}=\{\rho_{5}, \cdots ,
\rho_{8}\}.$ Let ${\mathcal E}_{i}$ be the operator associated with
each $G_{i}$. Then ${\mathcal E}_{2}$ is a unitary operator and
$${\mathcal E}_{i}={\mathcal E}_{2}^{i-1}$$
for $i\leq 4$ and ${\mathcal E}_{5}$ is self-adjoint.  Hence
$${\mathcal E}_{i}^{*}={\mathcal E}_{2}^{3(i-1)}.$$

For two local inverses $\tau_{1}$, $\tau_{2}$, if $\tau_{1}$ is
equivalent to $\tau_{2}$, then $\tau_{1}\circ\rho$ is equivalent to
$\tau_{2}\circ\rho$. Thus
$${\mathcal E}_{5}{\mathcal E}_{i}={\mathcal E}_{5}$$
for $i\leq 4$. So taking the adjoint of the above equality gives
$${\mathcal E}_{i}^{*}{\mathcal E}_{5}={\mathcal E}_{5}.$$
Therefore we have
$${\mathcal E}_{i}{\mathcal E}_{5}={\mathcal E}_{5},$$
to obtain that $\{{\mathcal E}_{1}, \cdots , {\mathcal E}_{4}\}$ is
in the center of ${\mathcal A}_{\phi}$ and hence ${\mathcal
A}_{\phi}$ is commutative.

If $G$ has $2$ elements, then the set of local inverses of $\phi$ is
divided into $G_{1}=\{\rho_{1}\}$, $G_{2}=\{\rho_{2}\}$,
$G_{3}=\{\rho_{3}, \rho_{4}\}$, $G_{4}=\{\rho_{5}, \rho_{6}\}$, and
$G_{5}=\{\rho_{7}, \rho_{8}\}.$ Let ${\mathcal E}_{i}$ be the
operator associated with $G_{i}$ for $i\leq 5$. We will show that
${\mathcal E}_{2}$ is in the center of ${\mathcal A}_{\phi}$. If this is true, we obtain that the dimension of the center of ${\mathcal
A}_{\phi}$ has dimension at least $2$ and hence ${\mathcal
A}_{\phi}$ is commutative.

Since  ${\mathcal E}_{i}^{*}$ is in $\{{\mathcal E}_{3}, {\mathcal
E}_{4} , {\mathcal E}_{5}\}$, we have that  either every operator in
$\{{\mathcal E}_{3}, {\mathcal E}_{4} , {\mathcal E}_{5}\}$ is
self-adjoint or only one of them is self-adjoint. Now we consider
two cases.

If  every operator in $\{{\mathcal E}_{3}, {\mathcal E}_{4} , {\mathcal
E}_{5}\}$ is self-adjoint, we observe that for two local inverses
$\tau_{1}$, $\tau_{2}$, if $\tau_{1}$ is equivalent to $\tau_{2}$,
then $
 \tau_{1}\circ\rho_{2}$ is equivalent to
$ \tau_{2}\circ\rho_{2}$ to get that there is a permutation $\sigma$
of $\{3,4, 5\}$ such that
$${\mathcal E}_{i}{\mathcal E}_{2}={\mathcal E}_{\sigma (i)}.$$
Taking the adjoint of the above equality gives
$${\mathcal E}_{2}{\mathcal E}_{i}={\mathcal E}_{\sigma (i)}.$$
Thus ${\mathcal E}_{2}$ is in the center of ${\mathcal A}_{\phi}$.

If only one of ${\mathcal E}_{3}, {\mathcal E}_{4} , {\mathcal
E}_{5}$ is self-adjoint, we may assume that ${\mathcal E}_{3}$ is
self-adjoint.  For two local inverses $\tau_{1}$, $\tau_{2}$, if
$\tau_{1}$ is equivalent to $\tau_{2}$, then
$\rho_{2}\circ\tau_{1}\circ\rho_{2}$ is equivalent to
$\rho_{2}\circ\tau_{2}\circ\rho_{2}$. Thus ${\mathcal
E}_{2}{\mathcal E}_{i}{\mathcal E}_{2}$ is in $\{{\mathcal E}_{3},
{\mathcal E}_{4} , {\mathcal E}_{5}\}$ for $3\leq i\leq 5.$ So there
is a permutation $\sigma $ of $\{3,4,5\}$ such that
$${\mathcal E}_{2}{\mathcal E}_{i}{\mathcal E}_{2}={\mathcal E}_{\sigma_{2} (i)},$$
and
$${\mathcal E}_{2}{\mathcal E}_{3}{\mathcal E}_{2}={\mathcal E}_{3}.$$
This implies that $\sigma_{2}(3)=3,$ and hence $\sigma_{2}$ is a
permutation of $\{4, 5\}$.
   Noting that for each $3\leq i\leq 5$,
$${\mathcal E}_{i}^*\in \{{\mathcal E}_{3}, {\mathcal E}_{4} , {\mathcal
E}_{5}\},$$ we see that there is the other permutation $\sigma_{1}$
of $\{3,4,5\}$ with order $2$ such that
 $${\mathcal E}_{i}^*={\mathcal
E}_{\sigma_{1} (i)}$$ and $\sigma_{1} (3)=3,$ which means that
$\sigma_{1}$ is also a permutation of $\{4, 5\}.$ This also gives
$$ {\mathcal E}_{2}{\mathcal E}_{\sigma_{1}(i)}{\mathcal E}_{2}={\mathcal E}_{\sigma_{2}(i)}^{*}
={\mathcal E}_{\sigma_{1}(\sigma_{2} (i))}.$$ Thus
$${\mathcal E}_{2}{\mathcal E}_{i}{\mathcal E}_{2}={\mathcal E}_{\sigma_{1}(\sigma_{2} (\sigma_{1}^{-1}(i)))}.$$
 Since $\sigma_{1}$ and $\sigma_{2}$ are permutations of $\{4, 5\}$,
they commute with each other. Thus we have
$${\mathcal E}_{2}{\mathcal E}_{i}{\mathcal E}_{2}={\mathcal E}_{i},$$
which implies
$${\mathcal E}_{2}{\mathcal E}_{i}={\mathcal E}_{i}{\mathcal E}_{2},$$
since ${\mathcal E}_{2}^{2}=I.$ This implies that ${\mathcal E}_{2}$
is in the center of ${\mathcal A}_{\phi}$.

Now we will show that the number of minimal reducing subspaces of $M_{\phi}$ equals  the number of connected components of the Riemann surface $\phi^{-1}\circ\phi$ over the unit disk.
Recall that a projection $E$ in a $C^{*}$-algebra is called minimal if the projection $E$ does not equal $0$ and the only subprojections of $E$ in the $C^{*}$-algebra are $0$ and $E$. Then for every minimal reducing subspace $\mathcal M$ of $M_{\phi}$, we can associate a minimal projection
$P_{\mathcal M}$ in the commutant $\{M_{\phi}\}^{\prime}$ where $P_{\mathcal M}$ is the orthogonal projection from the Bergman space $L^{2}_{a}$ onto $\mathcal M$.

Let ${\mathcal P}$ be the set of all minimal projections in the commutant $\{M_{\phi}\}^{\prime}.$  The Fuglede theorem \cite{Con}, \cite{Dou} gives that
 ${\mathcal P}$ is  the set of all minimal projections in the $C^{*}$-algebra $\{M_{\phi}\}^{\prime}\cap \{M_{\phi}^{*}\}^{\prime}.$
Let ${\mathcal A}$
be the $C^{*}$-algebra generated by those elements in ${\mathcal P}$. Hence $\mathcal A$ equals ${\mathcal A}_{\phi}$. As shown above, ${\mathcal A} $ is a commutative $C^{*}$-algebra and its dimension equals $q$.
Thus for any two minimal projections $P_{1}$ and $P_{2}$ in $\mathcal P$, $P_{1}$ commutes with $P_{2}$. Noting that both $P_{1}$ and $P_{2}$
are minimal, we have that $P_{1}P_{2}=0$ to get that $P_{1}$ is orthogonal to $P_{2}.$ Therefore $q$ equals the number of elements in ${\mathcal P}$.
Let ${\mathcal M}_{P}$ denote the range of the projection $P$. Thus ${\mathcal M}_{P}$ is a minimal reducing subspace of $M_{\phi}$.
So  the number of minimal reducing subspaces of $M_{\phi}$   equals $q$. This completes the proof.

\bigskip

{\bf Acknowledgement.} We thank Kunyu Guo and Hansong Huang for
their  comments on the earlier version of this paper.

\end{document}